\documentclass[12pt,a4paper]{amsart}
\usepackage{amsmath}
\usepackage{lscape}
\usepackage{graphics}

\newcommand{\mr}{\mathbb{R}} 

\begin{document}
\newcommand{\mc}{\mathbb{C}}
\newcommand{\mh}{\mathbb{H}}
\newcommand{\mo}{\mathbb{O}}
\newcommand{\mk}{\mathbb{K}}
\newcommand{\mj}{\mathbb{J}}
\newcommand{\mcs}{\tilde{\mathbb{C} }}
\newcommand{\mhs}{\tilde{\mathbb{H} }}
\newcommand{\mos}{\tilde{\mathbb{O} }}
\newcommand{\mks}{\tilde{\mathbb{K} }}
\newcommand{\sa}{\mathfrak{sa} }
\newcommand{\so}{\mathfrak{so} }
\newcommand{\su}{\mathfrak{su} }
\newcommand{\sq}{\mathfrak{sq} }
\newcommand{\g}{\mathfrak{g} }
\newcommand{\f}{\mathfrak{f} }
\newcommand{\p}{\mathfrak{p} }
\renewcommand{\aa}{\alpha}
\renewcommand{\gg}{\gamma}
\renewcommand{\d}{\delta }
\newcommand{\x}{\bar{x} }
\newcommand{\y}{\bar{y} }
\newcommand{\s}{\bar{s} }
\newcommand{\ds}{\dot{+}}
\newcommand{\Der}{\operatorname{Der}}
\newcommand{\Str}{\operatorname{Str}}
\newcommand{\ad}{\operatorname{ad}}
\newcommand{\Con}{\operatorname{Con}}
\newcommand{\Cl}{\operatorname{Cl}}
\renewcommand{\sp}{\mathfrak{sp} }
\renewcommand{\sl}{\mathfrak{sl}}
\newcommand{\tr}{\operatorname{tr}}
\newcommand{\diag}{\operatorname{diag}}
\newcommand{\Tri}{\operatorname{Tri}}
\renewcommand{\Im}{\operatorname{Im}}
\renewcommand{\Re}{\operatorname{Re}}
\renewcommand{\o}{\mathfrak{o} }
\renewcommand{\O}{\operatorname{O}}
\newcommand{\SO}{\operatorname{SO}}
\newcommand{\I}{\leavevmode\hbox{\rm\small1\kern-3.8pt\normalsize1}} 
\theoremstyle{plain} \newtheorem{prop}{Proposition}
\theoremstyle{plain} \newtheorem{theorem}{Theorem}
\theoremstyle{plain} \newtheorem{lemma}{Lemma}
\theoremstyle{definition} \newtheorem{defn}{Definition}
\today \\
\title{Magic squares of Lie Algebras}
\author{C.H. Barton}
\curraddr{Department of Mathematics\\ University of York\\ Heslington \\
York \\ YO10 5DD}
\email[C.H. Barton]{chb108@york.ac.uk}
\author{A. Sudbery}
\email[A. Sudbery]{as2@york.ac.uk}
\begin{abstract}
This paper is an investigation of the relation between Tit's magic
square of Lie algebras and certain Lie algebras of $3\times 3$ and
$6\times 6$ matrices with entries in alternative algebras. By
reformulating Tit's definition in terms of \emph{trialities} (a
generalisation of derivations), we give a systematic explanation of the
symmetry of the magic square. We show that when the columns of the magic
square are labelled by the real division algebras and the rows by their
split versions, then the rows can be interpreted as analogues of the
matrix Lie algebras $\su (3)$, $\sl (3)$ and $\sp (6)$ defined for each
division algebra. We also define another magic square based on $2\times
2$ and $4\times 4$ matrices and prove that it consists of various
orthogonal or (in the split case) pseudo-orthogonal Lie algebras.
\end{abstract}

\maketitle

\section{Introduction}
\label{sec:intro}
Semisimple Lie groups and Lie algebras are normally discussed in terms
of their root systems, which makes possible a unified treatment and
emphasises the common features of their underlying structures. However,
some classical investigations~\cite{Weyl} depend on particularly simple
matrix descriptions of Lie groups. This creates a
distinction between the classical groups (naturally enough) and the exceptional ones, which
is maintained in some more recent work (e.g.~\cite{Howe,Molev,Nazarov}).
This paper is motivated by the desire to give a similar matrix
description of the exceptional groups, thus assimilating them to the
classical groups, with a view to extending results like the Capelli
identities to the exceptional cases.

It has long been known~\cite{Molev} that most exceptional Lie algebras
are related to the exceptional Jordan algebra of $3\times 3$ hermitian
matrices with entries from the octonions, $\mo$. Here we show that this
relation yields descriptions of certain real forms of the complex Lie
algebras $F_{4}, E_{6}$ and $E_{7}$ which can be interpreted as
octonionic versions of, respectively, the Lie algebra of antihermitian
$3\times 3$ matrices, that of special linear $3\times 3$ matrices and
that of symplectic $6\times 6$ matrices. To be precise, we define for
each alternative algebra $\mk$ a Lie algebra $\sa (3,\mk )$ such that
$\sa (3,\mc )=\su(3)$ and $\sa(3,\mo)$ is the compact real form of
$F_{4}$; a Lie algebra $\sl(3,\mk)$ equal to $\sl (3,\mc)$ for $\mk=\mc$
and a non-compact real form of $E_{6}$ for $\mk=\mo$; and a Lie algebra
$\sp(6,\mk)$ such that  $\sp(6,\mc)$ is the set of $6\times 6$ complex
matrices $X$ satisfying $X^{\dag}J=-JX$, (where J is an antisymmetric
real $6\times 6$ matrix and $X^{\dag}$ denotes the hermitian conjugate
of $X$), and such that $\sp(6,\mo)$ is a non-compact real form of
$E_{7}$.

Our definitions can be adapted to yield Lie algebras $\sa(2,\mk), \sl
(2,\mk)$ and $\sp(4,\mk)$ reducing to $\su(2),\sl(2,\mc)$ and
$\sp(4,\mc)$ when $\mk=\mc$. These Lie algebras are isomorphic to
various pseudo-orthogonal algebras.

These constructions are all related to Tits's magic square of Lie
algebras~\cite{Tits62} and based on an unpublished suggestion of Ramond~\cite{Ramond}.
The magic square of Tits is a construction of a Lie algebra $L(\mj,\mk)$
for any Jordan algebra $\mj$ and alternative algebra $\mk$. If
$\mj=H_{3}(\mk_{1})$ is the Jordan algebra of $3\times 3$ matrices with
entries from an alternative algebra $\mk_{1}$ and if $\mk=\mk_{2}$ is
another alternative algebra, this yields a Lie algebra
$L_{3}(\mk_{1},\mk_{2})$ for any pair of alternative algebras. Taking
$\mk_{1}$ and $\mk_{2}$ to be real division algebras, we obtain a
$4\times 4$ square of compact Lie algebras which (magically) is
symmetric and contains the compact real forms of $F_{4}, E_{6}, E_{7}$
and $E_{8}$. We will show that if the division algebra $\mk_{2}$ is
replaced by its split form $\mks_{2}$, one obtains a non-symmetric
square of Lie algebras whose first three rows are the sets of matrix Lie
algebras described above:
\begin{align}
\label{eqn:iso2}
	L_{3}(\mk ,\mr )&=\sa (3,\mk) \notag \\
	L_{3}(\mk, \mcs )&=\sl (3,\mk) \\
	L_{3}(\mk ,\mhs )&=\sp (6,\mk). \notag
\end{align}
We will also describe magic squares of Lie algebras based on $2\times 2$
matrices, which have similar properties.

The organisation of the paper is as follows. In
Section~\ref{sec:notation} we establish notation and recall the
definitions of the various kinds of algebra with which we will be
concerned. In Section~\ref{sec:ms} we give Tits's definition of the Lie
algebras $L_{3}(\mk _{1},\mk _{2})$ and state the main properties of the
magic square; we also give the definition and properties of the $2\times
2$ magic square $L_{2}(\mk _{1},\mk _{2})$. Section ~\ref{sec:symmetry}
is concerned with the symmetry property of the $3\times 3$ magic square:
we reformulate the definition of $L_{3}(\mk _{1},\mk _{2})$, using
Ramond's concept of a \emph{triality algebra}, so as to make the
symmetry manifest. Section~\ref{sec:proofs2} contains proofs of the
properties of the $2\times 2$ magic square which were stated in
Section~\ref{sec:ms}, and Section~\ref{sec:proofs3} contains the
proofs of the corresponding properties of the $3\times 3$ magic square.

\section{Notation.}
\label{sec:notation}
We will use the notation $\ds$ to denote the direct sum of vector
spaces. This enables us to reserve the use of $\oplus$ to denote the
direct sum of Lie algebras, i.e. $A\oplus B$ implies that $[A,B]=0$.

An algebra $\mk$ (over $\mr$) with a non-degenerate quadratic form,
which we will denote by $x\mapsto \left| x\right| ^{2}$, satisfying
\begin{equation}
\label{eqn:qform}
	\left| xy\right| ^{2}=\left| x\right| ^{2}\left| y\right| ^{2} \quad x,y\in
\mk ,
\end{equation}
is known as a \textit{composition algebra}. We consider $\mr$ to be
embedded in $\mk$ as the set of scalar multiples of the identity
element, and denote by $\mk ^{\prime}$ the subspace of $\mk$ orthogonal
to $\mr$. It can then be shown~\cite{Kantorsolod} that $\mk=\mr \ds \mk^{\prime }$
 and we write $x=\Re x + \Im x$ with $\Re x \in \mr$ and $\Im x \in \mk
^{\prime}$. It can also be shown that the conjugation which fixes each
element of $\mr$ and multiplies every element of $\mk ^\prime $ by $-1$,
denoted $x\mapsto \overline{x}$, satisfies
\begin{equation}
	\overline{xy}=\overline{y}\, \overline{x}
\end{equation}
as well as
\begin{equation}
	x\overline{x}=\left| x\right| ^{2}.
\end{equation}
We use the notation $[x,y,z]$ for the associator
\begin{equation}
	[x,y,z] =(xy)z-x(yz).
\end{equation}
Any composition algebra $\mk$ satisfies the \textit{alternative law}, i.e. the
associator is an alternating function of $x,y$ and $z$. If $\left| x\right| ^{2}$
is positive definite then $\mk$ is a division algebra.
 
A division algebra is an algebra in which we have
\begin{equation*}
	xy=0 \Rightarrow  x=0  \text{ or }  y=0.
\end{equation*} The only such positive definite composition algebras are  $\mr ,\mc
,\mh $ and $\mo $
(Hurwitz's Theorem)~\cite{Schafer66} which we
denote in general by $\mk $. We denote their dimension by $\nu$, thus
$\nu = 1,2,4$ or $8$. These algebras are obtained from the
Cayley-Dickson process~\cite{Schafer66} and, using the same process with
different signs, \textit{split} forms
of these algebras can also be obtained. These are so called because in   $\mc ,\mh
$ and $\mo
$ we have
\begin{equation*}
	i^{2}+1=0
\end{equation*}
but in the split algebras $\mcs, \mhs $ and $\mos$ we have
\begin{equation*}
	i^{2}-1=(i+1)(i-1)=0
\end{equation*}
i.e. the equation can be \textit{split } for at least one of the imaginary
basis elements. Thus whilst the positive
definite algebras $\mr ,\mc ,\mh $ and $\mo $ are division algebras, the
split forms $\mcs, \mhs $ and $\mos$ are not.

Our notation for Lie algebras is that used in~\cite{Sudbery84}. We
use the notation $A^{\dag }$ for the hermitian conjugate of the
matrix $A$ with entries in $\mk$, defined in analogy to the complex case
by
\begin{equation*}
	(X^{\dag})_{ij}=\overline{X}_{ji}.
\end{equation*}
We use $\su (s,t)$ for the Lie algebra of the pseudo-unitary group, 
\begin{equation*}
	\su (s,t)=\{ A \in \mc ^{n\times n} :A^{\dag}G+GA=0\}
\end{equation*} 
where $G=\diag (-1,\dots ,-1,+1,\dots ,+1)$ with $s$ $-$ signs and $t$ $+$ signs;
$\sq (n)$ for the Lie algebra of antihermitian quaternionic matrices $A$,
\begin{equation*}
	\sq (n)=\{ A\in \mh ^{n\times n}: A^{\dag}=-A\} ;
\end{equation*}
and $\sp (2n,\mk )$ for the Lie algebra of the symplectic group of $2n\times
2n$ matrices with entries in $\mk$, i.e.
\begin{equation*}
	 \sp (2n,\mk )=\{ A\in \mk ^{2n\times 2n}:A^{\dag }J+JA=0\}
\end{equation*}
where $J=\begin{pmatrix} 0 & I_{n} \\ -I_{n} & 0 \\ \end{pmatrix}$.
We also have $\so (s,t)$, the Lie algebra of the pseudo-orthogonal
group $\SO (s,t)$, given by 
\begin{equation*}
	\so (s,t)=\{ A\in \mr ^{n\times n}:A^{T}G+GA=0\}
\end{equation*}
where $G$ is defined as before. We will also write $\O (V,q)$ for the
group of linear maps of the vector space $V$ preserving the
non-degenerate quadratic form $q$, $\SO (V,q)$ for its unimodular (or
special) subgroup, $\o (V,q)$ and $\so (V,q)$ for their Lie algebras. We
omit $q$ if it is understood from the context. Thus for any division
algebra we have $\SO (\mk )$ and $\so (\mk )$.

A Jordan algebra $\mj$ is defined to be a commutative algebra (over a
field $\mk$) in which all products satisfy the Jordan identity
\begin{equation*}
	(xy)x^{2}=x(yx^{2}).
\end{equation*}
Let $L_{n}(\mk )$ be the set of all $n\times n$ matrices with entries in
$\mk$, and let  $H_{n}(\mk )$ and $A_{n}(\mk )$ be the sets of all
hermitian and antihermitian matrices with entries in $\mk $
respectively. We denote by $H_{n}^{\prime}(\mk )$, $A_{n}^{\prime}(\mk
)$ and $L_{n}^{\prime}(\mk )$  the subspaces of traceless
matrices  of  $H_{n}(\mk )$, $A_{n}(\mk )$ and $L_{n}(\mk )$ respectively. 
We thus have $L_{n}(\mk )=H_{n}(\mk )\ds A_{n}(\mk )$ and $L_{n}^{\prime }(\mk
)=H_{n}^{\prime }(\mk )\ds A_{n}^{\prime }(\mk )$. We will use the
fact that $H_{n}(\mk )$ is a Jordan algebra for $\mk = \mr ,\mc ,\mh$
for all $n$ and for 
$\mk=\mo$ when $n=2,3$~\cite{Sudbery84}, with the Jordan product as the anticommutator
\begin{equation*}
	X\cdot Y=XY+YX.
\end{equation*}
This is a commutative but non-associative product. 

The derivation algebra, $\Der A$, of any algebra $A$ is defined as
\begin{equation}
\label{eqn:Der}
	\Der A = \{D \mid  D(xy)=D(x)y + xD(y) \}
\end{equation}
for $x,y \in A$. 
The derivation algebras of the four positive definite composition
algebras are as follows:
\begin{align*}
	\Der \mr &= \Der \mc =0, \\
	\Der \mh &= C(\mh ^{\prime})=\{ C_{a} \mid a \in \mh ^{\prime}
\} \text{ where } C_{a}(q)=aq-qa. 
\end{align*}
$\Der \mo $ is an exceptional Lie algebra of type $G_{2}$.

The structure algebra $\Str A$ of any algebra $A$ is
defined to be the Lie algebra generated by left and right multiplication
maps $L_a$ and $R_{a}$ for $a\in A$. For Jordan algebras this can be shown to be~\cite{Schafer66}
\begin{equation}
\label{eqn:Str}
	\Str  \mathbb{J} = \Der \mathbb{J} \ds L(\mathbb{J})
\end{equation}
where $L(\mathbb{J} )$ is the set of all $L_a$ with $a\in \mathbb{J}$.
The algebra denoted by $\Str ^{\prime} \mathbb{J}$ 
is the structure algebra with its centre factored out. 
We also require another Lie algebra
 associated with a Jordan algebra, namely the conformal algebra as
constructed by Kantor (1973) and Koecher (1967). The underlying vector
space of this is
\begin{equation}
\label{eqn:Con}
	\Con  \mathbb{J} = \Str  \mathbb{J} \ds \mathbb{J} ^{2}.
\end{equation}
We will specify the Lie brackets of these algebras at a later stage.

\section{Magic Squares: Summary of Results}
\label{sec:ms}
\subsection{$3\times 3$ Matrices}
\label{sec:ms3}
Let $\mk$ be a real composition algebra and $\mj$ a real Jordan algebra, with $\mk
^{\prime}$ and $\mj ^{\prime}$ the quotients of the algebras by the
subspaces of scalar multiples of the identity. Define a vector space
\begin{equation}
	M(\mathbb{J} ,\mk )=\Der \mathbb{J} \ds (\mathbb{J} ^{\prime} \otimes
\mk ^{\prime })\ds \Der \mk.
\end{equation}
Then define 
\begin{equation*}
	L_3(\mk _1,\mk _2) = M(H_{3}(\mk _1) ,\mk _2).
\end{equation*}
Explicitly this is the vector space
\begin{equation}
\label{eqn:L3a}
	L_3(\mk _1,\mk _2) = \Der H_3(\mk _1)\ds H^\prime _3(\mk _1)\otimes
\mk ^\prime _2 \ds \Der \mk _2
\end{equation}
which is a Lie algebra with Lie subalgebras $\Der H_{3}(\mk _{1})$ and
$\Der \mk_{2}$ when taken with the brackets \label{brackets}
\begin{align}
\label{eqn:Tits}
	[D,A\otimes x] &= D(A)\otimes x   \notag  \\
 	[E,A\otimes x] &= A\otimes E(x)    \notag \\
	[D,E] &= 0 \\
	[A\otimes x,B\otimes y] &= \tfrac{1}{6} \langle A,B\rangle
D_{x,y}+(A\ast B)\otimes \tfrac{1}{2}[x,y]-\langle x,y\rangle [L_{A},L_{B}] \notag
\end{align}
with $D \in \Der H_3(\mathbb{K}_{1});\,  A,B\in H^{\prime }_3(\mathbb{K}
_1);\, x,y\in \mathbb{K} ^{\prime }_2$ and $E\in \Der \mathbb{K} _2$. 
These brackets are obtained from Schafer's description of the Tits
construction~\cite{Schafer66}. They require some explanation.  $\langle A,B\rangle
$ and $(x,y)$ denote the symmetric bilinear forms on $H_{3}(\mk_{1})$ and
$\mk_{2}$ respectively, given by
\begin{align*}
	\langle A,B\rangle &=\Re (\tr (A\cdot B))=2\Re (\tr (AB)) \\
	\langle x,y\rangle &=\tfrac{1}{2}  (\left| x+y\right| ^{2}-\left| x\right| ^{2}-\left| y\right|
^{2}) = \Re (x\overline{y}).
\end{align*}
The derivation $D_{x,y}$ is defined as
\begin{equation}
\label{eqn:Ddef}
	D_{x,y}=[L_{x},L_{y}] +[L_{x},R_{y}] +[R_{x},R_{y}] \quad \in
\Der \mk _{2}.
\end{equation}
For future reference we note that 
\begin{equation}
\label{eqn:Deq}
	D_{x,y}z=[[x,y],z]-3[x,y,z]
\end{equation}
which shows that $D_{x,y}=-D_{y,x}$.
Finally $(A\ast B)$ is the traceless part of the Jordan product of $A$
and $B$,
\begin{equation}
	A\ast B = A\cdot B -\tfrac{1}{3}\tr (A\cdot B).
\end{equation}

Tits~\cite{Tits62} (see also~\cite{Freudenthal65,Schafer66}) showed that this gives a unified
construction leading to the so-called magic square of Lie algebras of
$3\times 3$ matrices whose complexifications are
\begin{center}
\vspace{0.5cm} \begin{tabular}{|c||c|c|c|c|}
\hline
   & $\mr$ & $\mc$  & $\mh$ & $\mo$ \\
 \hline \hline
  $\mr$ & $A_1$ & $A_2$ & $C_3$ & $F_4$ \\
  \hline
   $\mc$ & $A_2$ & $A_2 \oplus A_2$ & $A_5$ & $E_6$ \\
   \hline
    $\mh$ & $C_3$ & $A_5$ & $B_6$ & $E_7$ \\
    \hline
    $\mo$ & $F_4$ & $E_6$ & $E_7$ & $E_8$ \\
	\hline
    \end{tabular} \vspace{0.5cm}
    \end{center} 
The striking properties of this square are (a) its symmetry and (b) the
fact that four of the five exceptional Lie algebras occur in its last
row. The explanation of the symmetry property is the subject of
section~\ref{sec:symmetry}.
The fifth exceptional Lie algebra, $G_{2}$, can be included by adding an extra row
corresponding to the Jordan algebra
$\mr$. 

In~\cite{Sudbery84} it is asserted without proof that we can write this in a
slightly different form and
that we can include the isomorphisms listed in~\eqref{eqn:iso2}. 
This involves a different set of real forms obtained by taking  the
split composition algebras $\mr ,\mcs, \mhs ,$
and $\mos$ rather than $\mr ,\mc ,\mh,$ and $\mo$ as the second algebra. Thus the
split magic square for three by
three matrices looks like
\begin{center}
\vspace{0.5cm} \begin{tabular}{|c||c|c|c|c|}
\hline
     & $\mr$ & $\mc$  & $\mh$ & $\mo$ \\
 \hline \hline
  $\Der H_3(\mk )\cong L_3(\mk ,\mr)$ & $\so(3)$ & $\su(3)$ & $\sq(3)$ & $F_{4,1}$ \\
  \hline
   $\Str ^\prime  H_3(\mk )\cong L_3(\mk ,\mcs ) $ & $\sl(3,\mr )$ &
$\sl(3,\mc )$ & $\sl(3,\mh )$ & $E_{6,1}$ \\
   \hline
    $\Con  H_3(\mk )\cong L_3(\mk ,\mhs ) $ & $\sp(6,\mr )$ & $\su(3,3)$ &
$\sp(6,\mh )$ & $E_{7,1}$ \\
    \hline
    $L_3(\mk ,\mos) $ & $F_{4,2}$ & $E_{6,2}$ & $E_{7,2}$ & $E_{8,1}$ \\
	\hline
    \end{tabular} \vspace{0.5cm}
    \end{center}
where the notation $_{,1}$ and $_{,2}$ (in the style
of~\cite{Freudenthal65}) is used to distinguish between
different real forms of the exceptional Lie algebras in the last row and
column. These are identified by their maximal compact subalgebras
as follows:
\begin{center}
\label{maximal}
\vspace{0.5cm} \begin{tabular}{|c|c|}
\hline
	Exceptional Lie Algebra & Maximal Compact Subalgebra \\
\hline

	$E_{6,1}$ & $F_4$ \\
\hline
	$E_{7,1}$ & $E_{6,1} \oplus \so (2)$ \\
\hline
	$E_{8,1}$ & $E_{7,1} \oplus \so (3)$ \\
\hline 

	$E_{6,2}$ & $\sq (3) \oplus \so (3)$ \\
\hline
	$E_{7,2}$ & $\su (6) \oplus \so (3)$ \\
\hline
	$E_{8,2}$ & $\so (12) \oplus \so (3)$ \\
\hline 

\end{tabular} \vspace{0.5cm}.
\end{center}

\subsection{$2\times 2$ Matrices}

The Tits construction can also be adapted for $2\times 2$ matrix
algebras. In this case we take the vector space to be
\begin{equation}
	L_2(\mk _{1},\mk _{2})=\Der H_{2}(\mk _{1})\ds H^{\prime }_{2}(\mk
_{1})\otimes \mk ^{\prime }_{2}\ds \so (\mk ^{\prime }_{2})
\end{equation}
which is again a Lie algebra when taken with the brackets
\begin{align}
\label{eqn:2bracs}
	[D,A\otimes x] &=D(A)\otimes x   \\
	[E,A\otimes x] &=A\otimes E(x) \notag \\
	[D,E] &= 0 \notag\\
	[A\otimes x,B\otimes y] &=\tfrac{1}{4} \langle A,B\rangle
D_{x,y}-\langle x,y\rangle [R_A,R_B] \notag
\end{align}
where the symbols used in this set of brackets are defined in the same
way as the ones used in the $3\times 3$ case. We note that
$D_{x,y}=2s_{x,y}$, where $s_{x,y}$ is the element of $\so (\mk ^{\prime
}_{2})$ that maps $x$ to $y$ and $y$ to $\pm x$, depending on the metric
of $\mk ^{\prime }_{2}$ i.e.
\begin{equation}
\label{eqn:Sdef}
	S_{x,y}(z)=\langle x,z\rangle y - \langle y,z \rangle x
\end{equation}

If $\mk_{1},\mk_{2}$
are division algebras then  this gives the compact magic square for $2\times 2$ matrix
algebras
\begin{equation*}
	L_{2}(\mk_{1},\mk_{2})=\so(\nu_{1}+\nu_{2}).
\end{equation*}
If $\mk_{2}$ is one of the split composition algebras $\mcs,\mhs$ or
$\mos$ this becomes
\begin{equation*}
	L_{2}(\mk_{1},\mk_{2})=\so(\nu_{1}+\tfrac{1}{2}
\nu_{2},\tfrac{1}{2} \nu_{2}).
\end{equation*}
giving the magic square
\begin{center}
\vspace{0.5cm} \begin{tabular}{|c||c|c|c|c|} \hline
	& $\mr$ & $\mc$ & $\mh$ & $\mo$ \\ \hline \hline
$L_2(\mk ,\mr )$ & $\so (2)$ & $\so (3)$ & $\so (5)$ & $\so (9)$ \\ \hline
$L_2(\mk ,\mcs )$ & $\so (2,1)$ & $\so (3,1)$ & $\so (5,1)$ & $\so (9,1)$ \\ \hline
$L_2(\mk ,\mhs )$ & $\so (3,2)$ & $\so (4,2)$ & $\so (6,2)$ & $\so (10,2)$ \\ \hline
$L_2(\mk ,\mos )$ & $\so (5,4)$ & $\so (6,4)$ & $\so (8,4)$ & $\so (12,4)$ \\ \hline
\end{tabular} \vspace{0.5cm}.
\end{center}
As in the $3\times 3$ case, these Lie algebras can be identified with
certain types of $2\times 2$ matrix algebras
\begin{center}
\vspace{0.5cm} \begin{tabular}{|c||c|c|c|c|} \hline
	& $\mr$ & $\mc$ & $\mh$ & $\mo$ \\ \hline \hline
$\Der H_2(\mk )\cong L_2(\mk ,\mr )$ & $\so(2)$ & $\su(2)$ & $\sq(2)$ & $\so(9)$ \\ \hline
$\Str  H_2(\mk )\cong L_2(\mk ,\mcs )$ & $\sl(2,\mr )$ & $\sl(2,\mc )$ &
$\sl(2,\mh )$ & $\sl(2,\mo )$ \\ \hline
$\Con  H_2(\mk )\cong L_2(\mk ,\mhs )$ & $\sp(4,\mr )$ & $\su(2,2)$ &
$\sp(4,\mh )$ & $\sp(4,\mo )$ \\ \hline
$L_2(\mk ,\mos )$ & $\so(5,4)$ & $\so(6,4)$ & $\so(8,4)$ & $\so(12,4)$ \\ \hline
\end{tabular} \vspace{0.5cm}.
\end{center}
Again this extends the concepts of the Lie algebras
$\sa(2,\mk),\sl(2,\mk)$ and $\sp(2,\mk)$ to $\mk=\mh$ and $\mo$. Note
that $\su (2,2) \cong \sp (4,\mc )$.

\section{Symmetry Property of the $3\times 3$ magic square.}
\label{sec:symmetry}
In this section we will rearrange the definition
\begin{equation}
\label{eqn:L3}
	L_{3}(\mk _{1},\mk _{2})=\Der H_{3}(\mk_{1}) \ds H_{3}^{\prime }
(\mk _{1}) \otimes \mk _{2}^{\prime } \ds \Der \mk _{2}
\end{equation}
so as to make explicit the symmetry between $\mk _{1}$ and $\mk _{2}$.
We need a new Lie algebra associated with any $\mk$, defined as follows:
\begin{defn}
Let $\mk $ be a composition algebra over $\mr$. The \emph{triality
algebra} of $\mk $ is
\begin{equation}
	\Tri \mk = \{ (A,B,C) \in \so (\mk )^{3}| A(xy)=x(By)+(Cx)y,
\forall x,y\in \mk \}.
\end{equation}
\end{defn}
It is easy to verify that $\Tri \mk$ is a Lie algebra with brackets
defined componentwise, i.e. it is a Lie subalgebra of $\so (\mk ) \oplus \so (\mk ) \oplus \so (\mk
)$.

\begin{lemma}
The triality algebras of the four positive-definite composition
algebras can be identified as follows:
\begin{align*}
	\Tri \mr &=0 \\
	\Tri \mc &\cong \mr ^{2} \\
	\Tri \mh &\cong \so (3) \oplus \so(3) \oplus \so (3) \\
	\Tri \mo &\cong \so (8)
\end{align*}

\begin{proof}
$\Tri \mr =0$ because $\so (\mr)=0$. For $\mc $ we can identify $\so
(\mc )$ with the set of multiplication maps $z\mapsto hz$ with $h$ pure
imaginary, which is isomorphic to $\mr$ as a Lie algebra. Then $\Tri
\mc$ is the subspace of the abelian Lie algebra $\mr ^{3}$ given by
\begin{equation*}
	\Tri \mc = \{ (u,v,w): u=v+w \}
\end{equation*}
which is two dimensional.

Antisymmetric linear maps $A:\mh \to \mh $ are all of the form
$A=L_{a_{1}}+R_{a_{2}}$ with $a_{1},a_{2} \in \mh ^{\prime}$. These
are all independent (this is a reflection of the Lie algebra isomorphism
$\so (4) \cong \so (3) \oplus \so (3)$). Hence the condition for $(A,B,C)
\in \Tri \mh$ is of the form
\begin{equation*}
	a_{1}xy+xya_{2}=c_{1}xy+x(c_{2}+b_{1})y+xyb_{2}, \qquad \forall
x,y \in \mh
\end{equation*}

Taking $y=1$ and using the independence of the left and right
multiplication maps gives
\begin{equation*}
	a_{1}=c_{1} \quad \text{and} \quad a_{2}=c_{2}+b_{1}+b_{2}.
\end{equation*}
Taking $x=1$ gives
\begin{equation*}
a_{1}=c_{1}+c_{2}+b_{1} \quad \text{and} \quad a_{2}=b_{2}.
\end{equation*}
Hence $c_{2}+b_{1}=0$ and we have
\begin{equation*}
	A=L_{a_{1}}+R_{a_{2}}, \quad B=L_{b_{1}}+R_{a_{2}}, \quad
C=L_{a_{1}}-R_{b_{1}}.
\end{equation*}
Thus $\Tri \mh \cong \mh ^{\prime \, 3} \cong \so (3) \oplus \so (3)
\oplus \so (3)$.

Finally, the infinitesimal version of the principle of
triality~\cite{Porteous} asserts that for each $A\in \so (8)$ there are
unique $B,C \in \so (8)$ such that
\begin{equation*}
	A(xy)=x(By)+(Cx)y \qquad \forall x,y,\in \mo .
\end{equation*}
This establishes an isomorphism between $\Tri \mo $ and $\so (8)$. 
\end{proof}
\end{lemma}

We will now describe the chain of inclusions
\begin{equation}
	\Der \mk \subset \Tri \mk \subset \Der H_{3}(\mk )
\end{equation}
in a unified way, valid for any composition algebra $\mk$. We will use a
multiple notation to describe multiple direct sums, writing
\begin{equation*}
	nV = \underbrace{V \ds V \ds \dots \ds V}_{n}
\end{equation*}
rather than $V^{n}$ (which might suggest $V\otimes \dots \otimes V$).
\begin{lemma}
\label{lemma:pig}
For any composition algebra $\mk$,
\begin{equation*}
	\Tri \mk = \Der \mk \ds 2\mk ^{\prime}
\end{equation*}
in which $\Der \mk $ is a Lie subalgebra,
\begin{align*}
	[D,(a,b)] &= (Da, Db) \in 2\mk ^{\prime} \\
	[(a,0),(b,0)] &= \tfrac{2}{3}D_{a,b}+
\left( \tfrac{1}{3}[a,b],-\tfrac{2}{3}[a,b]\right), \\
	[(a,0),(0,b)] &= \tfrac{1}{3}D_{a,b}-\left( \tfrac{1}{3}[a,b],\tfrac{1}{3}[a,b]\right),
\\
	[(0,a),(0,b)] &= \tfrac{2}{3}D_{a,b}+\left(
-\tfrac{2}{3}[a,b],\tfrac{1}{3}[a,b]\right).
\end{align*}
\begin{proof}
Define $T:\Der \mk \ds 2\mk ^{\prime}\to \Tri \mk $ by
\begin{multline}
\label{eqn:Teq}
	T(D,a,b)= \\
	(D+L_{a}-R_{b},\, D-L_{a}-L_{b}-R_{b},\, D+L_{a}+R_{a}+R_{b}).
\end{multline}
This belongs to $\Tri \mk $ as a consequence of the alternative law. The
map $T$ is injective, for $T(D,a,b)=0$ implies
\begin{equation*}
	2L_{a}+L_{b}=0 \quad \text{and} \quad R_{a}+2R_{b}=0
\end{equation*}
(subtracting the first component from the second and third in turn), so
$2a+b=a+2b=0$ and hence $a+b=0$, which implies $D=0$. To show that $T$
is surjective, suppose $(A,B,C)\in \Tri \mk$ and define $a,b \in \mk
^{\prime}$ by
\begin{align}
	&B(1)= -a-2b, \notag \\
	&C(1)= 2a+b. \notag 
\end{align}
Let
\begin{equation}
\label{eqn:D}
	D = A-L_{a}+R_{b}.
\end{equation}
Since $(A,B,C) \in \Tri \mk$ we have
\begin{equation*}
	B(x)=1.B(x)=A(1.x)-C(1)x.
\end{equation*}
Thus
\begin{align}
\label{eqn:B}
	B&= A-2L_{a}-L_{b}  \notag\\
	&= D-L_{a}-L_{b}-R_{b} 
\end{align}
and
\begin{equation}
\label{eqn:C}
	C=D+L_{a}+R_{a}+R_{b} 
\end{equation}
Now 
\begin{align}
	D(xy) &= A(xy)-a(xy)+(xy)b \notag\\
	&= x(By)+(Cx)y-a(xy)+(xy)b \notag\\
	&=(Dx)y+x(Dy)
\end{align}
by equations~(\ref{eqn:B},~\ref{eqn:C}) and the alternative law. Hence $D$ is a derivation and
$(A,B,C)=T(D,a,b)$.

The first Lie bracket stated above, i.e.
\begin{equation*}
	[T(D,0,0),\, T(0,a,b)]=T(0,Da,Db),
\end{equation*}
follows from
\begin{equation}
	[D,L_{a}]=L_{Da} \quad \text{and} \quad [D,R_{a}]=R_{Da}.
\end{equation}
The other brackets follow from the commutators
\begin{align*}
	[L_{x},L_{y}]&=
\tfrac{2}{3}D_{x,y}+\tfrac{1}{3}L_{[x,y]}+\tfrac{2}{3}R_{[x,y]} \\
	[L_{x},R_{y}]&=
-\tfrac{1}{3}D_{x,y}+\tfrac{1}{3}L_{[x,y]}-\tfrac{1}{3}R_{[x,y]} \\
	[R_{x},R_{y}]&=
\tfrac{2}{3}D_{x,y}-\tfrac{2}{3}L_{[x,y]}-\tfrac{1}{3}R_{[x,y]}
\end{align*}
which can be calculated using equation~(\ref{eqn:Deq}).
\end{proof}
\end{lemma}

Any two elements $x,y\in \mk $ are associated with a derivation
$D_{x,y}\in \Der \mk $ and also with an antisymmetric map $S_{x,y}\in
\so (\mk)$, the generator of rotations in the plane of $x$ and $y$ given
by~(\ref{eqn:Sdef}).
There is also an element of $\Tri \mk$ associated with $x$ and $y$:
\begin{lemma}
\label{lemma:3}
For any $x,y\in \mk$, let
\begin{equation*}
	T_{x,y}=(4S_{x,y},\, 
R_{y}R_{\bar{x}}-R_{x}R_{\bar{y}},\, L_{y}L_{\bar{x}}-L_{x}L_{\bar{y}}).
\end{equation*}
Then $T_{x,y}\in \Tri \mk $.
\begin{proof}
We can write the action of $S_{x,y}$ as
\begin{align}
\label{eqn:cat}
	2S_{x,y}z &= (x\bar{z}+z\bar{x})y-x(\bar{z}y+\bar{y}z) \\
	&= -[x,y,z]+z(\bar{x}y)-(x\bar{y})z
\end{align}
using the alternative law and the relation $[x,y,\bar{z}]=-[x,y,z]$.
Since $\Re (\bar{x}y)=\Re (x\bar{y})$, we can write the last two terms as
\begin{align}
\label{eqn:dog}
	z(\bar{x}y)-(x\bar{y})z &= z\Im (\bar{x}y)-\Im (x\bar{y})z \\
	&=\tfrac{1}{2}z(\bar{x}y-\bar{y}x)-\tfrac{1}{2}(x\bar{y}-y\bar{x})z.
\end{align}
Now, by equation~(\ref{eqn:Deq}), we have
\begin{equation}
	S_{x,y}=\tfrac{1}{6}D_{x,y}+L_{a}-R_{b} 
\end{equation}
with
\begin{align*}
	a &= -\tfrac{1}{6}[x,y]-\tfrac{1}{4}(x\bar{y}-y\bar{x}) \, \in \mk
^{\prime} \\
	b &= -\tfrac{1}{6}[x,y]-\tfrac{1}{4}(\bar{x}y-\bar{y}x) \, \in \mk
^{\prime }.
\end{align*}
Hence, by equation~(\ref{eqn:Teq}), there is an element $(A,B,C) \in \Tri \mk $
with $A=S_{x,y}$ and
\begin{align*}
	B &= \tfrac{1}{6}D_{x,y}-L_{a}-L_{b}-R_{b}= S_{x,y}-L_{2a+b}, \\
	C &= \tfrac{1}{6}D_{x,y}+L_{a}+R_{a}+R_{b}= S_{x,y}+R_{a+2b},
\end{align*}
Writing $[x,y]=-\tfrac{1}{2}([\bar{x},y]+[x,\bar{y}])$ gives
\begin{align*}
	a+2b&=\tfrac{1}{4}(\bar{y}x-\bar{x}y) \\
	2a+b &=\tfrac{1}{4}(y\bar{x}-x\bar{y})
\end{align*}
so equations~(\ref{eqn:cat}) and~(\ref{eqn:dog}) give
\begin{equation}
	S_{x,y} = \tfrac{1}{2}Q_{x,y}-R_{a+2b}+L_{2a+b}
\end{equation}
where $Q_{x,y}z=-[x,y,z]$. Hence
\begin{align*}
	Cz &= -\tfrac{1}{2}[x,y,z]+\tfrac{1}{4}(y\bar{x}-x\bar{y})z \\
	&=\tfrac{1}{4}y(\bar{x}z)-\tfrac{1}{2}x(\bar{y}z)
\end{align*}
i.e.
\begin{equation*}
	C=\tfrac{1}{4}(L_{y}L_{\bar{x}}-L_{x}L_{\bar{y}}) 
\end{equation*}
and similarly
\begin{equation*}
	B=\tfrac{1}{4}(R_{y}R_{\bar{x}}-R_{x}R_{\bar{y}}).
\end{equation*}
Thus $T=(4S_{x,y},4C,4B)$ is an element of $\Tri \mk $.
\end{proof}
\end{lemma}
Note that if $x,y \in \mk ^{\prime}$, so that $\bar{x}=-x$ and
$\bar{y}=-y$, then $a=b=\tfrac{1}{12}[x,y]$ and so
\begin{multline}
\label{eqn:Teq2}
	T_{x,y}=( \tfrac{2}{3}D_{x,y}+\tfrac{1}{3}L_{[x,y]}-\tfrac{1}{3}R_{[x,y]}, 
 \tfrac{2}{3}D_{x,y}+\tfrac{1}{3}L_{[x,y]}+\tfrac{2}{3}R_{[x,y]},
 \\
 \tfrac{2}{3}D_{x,y}-\tfrac{2}{3}L_{[x,y]}-\tfrac{1}{3}R_{[x,y]} ) .
\end{multline}

The element $T_{x,y}$ will be needed to describe $\Tri \mk$ as a Lie
subalgebra of $\Der H_{3}(\mk )$. We will also need an automorphism of
$\Tri \mk$ defined as follows. For any linear map $A:\mk \to \mk$, let
$\overline{A}=KAK$, where $K:\mk \to \mk$ is the conjugation $x\mapsto
\bar{x}$ in $\mk$, i.e.
\begin{equation*}
	\overline{A}(x) = \overline{A(\bar{x})}.
\end{equation*}
Then $\overline{\overline{A}}=A$ and $\overline{AB}=\overline{A}\,
\overline{B}$. Note also that
\begin{equation*}
	\overline{L}_{x}=R_{\bar{x}}
\end{equation*}
and $D=\overline{D}$ if $D \in \so (\mk ^{\prime})$, in particular if D
is a derivation of $\mk$. 

\begin{lemma}
\label{lemma:splodge}
Given $T=(A,B,C) \in \Tri \mk$, let 
\begin{equation*}
	\theta (T) = (\overline{B},C,\overline{A}).
\end{equation*}
Then $\theta (T) \in \Tri \mk$ and $\theta$ is a Lie algebra
automorphism.
\begin{proof}
By Lemma~\ref{lemma:pig}, $T=T(D,a,b)$ for some $D\in \Der \mk$ and $a,b
\in \mk ^{\prime}$. Then
\begin{align*}
	A&= D+L_{a}-R_{b} \\
	B&= D-L_{a}-L_{b}-R_{b} \\
	C&= D+L_{a}+R_{a}+R_{b}.
\end{align*}
It follows that
\begin{equation*}
	\overline{B}=D+R_{a}+R_{b}+L_{b}=D+L_{a^{\prime}}-R_{b^{\prime
}}
\end{equation*}
which is the first component of
$T^{\prime}=(A^{\prime},B^{\prime},C^{\prime})\in \Tri \mk $, where
\begin{align*}
	B^{\prime} &=D-L_{a^{\prime}}-L_{b^{\prime}}-R_{b^{\prime}} \\
	&= D-L_{b}+L_{a+b}+R_{a+b}=C \\
	C^{\prime} &= D+L_{a^{\prime}}+R_{a^{\prime}}+R_{b^{\prime}} \\
	&=D+L_{b}+R_{b}-R_{a+b}=\overline{A},
\end{align*}
i.e. $T^{\prime}=(\overline{B},C,\overline{A})=\theta(T)$. It is clear
that $\theta$ is a Lie algebra automorphism.
\end{proof}
\end{lemma}
Given $T=(A,B,C)\in \Tri \mk $, it is convenient to define
$(T_{1},T_{2},T_{3})=(A,\overline{B},\overline{C})$. Then
$\theta (T)_{i} = T_{\sigma (i)}$ where $\sigma \in \mathcal{S}_{3}$ is
the cyclic permutation
\begin{equation*}
	\sigma (1)=2, \quad \sigma (2)=3, \quad \sigma (3)=1.
\end{equation*}

Write $a_{1}=-a-b$, $a_{2}=a$, $a_{3}=b$. Then the triality obtained
from $(a,b)$ can be written in the symmetric form
$T(0,a,b)=(T_{1},\overline{T}_{2},\overline{T}_{3})$ where
\begin{align*}
	T_{1}&=L_{a_{2}}-R_{a_{3}} \\
	T_{2}&=L_{a_{3}}-R_{a_{1}} \\
	T_{3}&=L_{a_{1}}-R_{a_{2}}
\end{align*}
i.e. $T_{i}=L_{a_{j}}-R_{a_{k}}$ where $(i,j,k)$ is a cyclic permutation
of $(1,2,3)$.

\begin{theorem}
\label{theorem:banana}
For any composition algebra $\mk$,
\begin{equation*}
	\Der H_{3}(\mk)=\Tri \mk \ds 3\mk
\end{equation*}
in which $\Tri \mk $ is a Lie subalgebra, and the brackets in $[ \Tri
\mk ,3\mk ]$ are
\begin{equation}
\label{eqn:1}
	[T,F_{i}(x)]=F_{i}(T_{i}x) \quad \in 3\mk ,
\end{equation}
if $T=(T_{1},\overline{T}_{2},\overline{T}_{3}) \in \Tri \mk$ and
$F_{1}(x)+F_{2}(y)+F_{3}(z)=(x,y,z) \in 3\mk $; and the brackets in
$[\Tri \mk ,\Tri \mk ]$ are given by 
\begin{equation}
\label{eqn:2}
	[F_{i}(x),F_{j}(y)]=F_{k}(\bar{y}\bar{x}) \in 3\mk,
\end{equation}
if $x,y \in \mk$ and $(i,j,k)$ is a cyclic permutation of $(1,2,3)$; and
\begin{equation}
\label{eqn:3}
	[F_{i}(x),F_{i}(y)]= \theta ^{1-i}(T_{x,y}) \in \Tri \mk .
\end{equation}
\begin{proof}
Define elements $e_{i},P_{i}(x)$ of $H_{3}(\mk)$ (where $i=1,2,3;\,  x \in
\mk $) by the equation
\begin{equation}
\label{eqn:Jdef}
\begin{pmatrix} \alpha & z & \bar{y} \\ \bar{z} & \beta & x \\ y &
\bar{x} & \gamma \end{pmatrix} = \alpha e_{1}+\beta e_{2}+\gamma e_{3} +
P_{1}(x)+P_{2}(y)+P_{3}(z) 
\end{equation}
for $\alpha , \beta , \gamma \in \mr ; \, x,y,z \in \mk $. Then the Jordan
product in $H_{3}(\mk)$ is given by
\begin{subequations}
\begin{align}
	e_{i}\cdot e_{j}&=2\delta _{ij} e_{i}  \label{eqn:J1} \\
	e_{i}\cdot P_{j}(x) &= (1-\delta _{ij})P_{j}(x) \label{eqn:J2}
\\
	P_{i}(x)\cdot P_{i}(y) &= 2(x,y)(e_{j}+e_{k}) \label{eqn:J3} \\
	P_{i}(x)\cdot P_{j}(y) &= P_{k}(\bar{y}\, \bar{x})
\label{eqn:J4}
\end{align}
\end{subequations}
where in each of the last two equations $(i,j,k)$ is a cyclic
permutation of $(1,2,3)$.

Now let $D:H_{3}(\mk)\to H_{3}(\mk)$ be a derivation of this algebra.
First suppose that
\begin{equation*}
	De_{i}=0, \quad i=1,2,3.
\end{equation*}
Then
\begin{align*}
	e_{i}\cdot DP_{i}(x) &=0 \\
	e_{i}\cdot DP_{j}(x) &= DP_{j}(x) \quad \text{if } i \neq j
\end{align*}
Thus $DP_{j}(x)$ is an eigenvector of each of the multiplication
operators $L_{e_{i}}$, with eigenvalue $0$ if $i=j$ and $1$ if $i\neq
j$. It follows that 
\begin{equation}
	DP_{j}(x)=P_{j}(T_{j}x)
\end{equation}
for some $T_{j}:\mk \to \mk $. Now
\begin{equation*}
	DP_{j}(x)\cdot P_{j}(y) +P_{j}(x)\cdot DP_{j}(y)=0
\end{equation*}
gives $T_{j} \in \so (\mk )$; and the derivation property of $D$ applied
to equation~(\ref{eqn:J4}) gives
\begin{equation*}
	T_{k}(\bar{y}\bar{x})=\bar{y}(\overline{T_{i}x})+(\overline{T_{j}y})\bar{x} 
\end{equation*}
i.e. $(T_{k},\overline{T_{i}},\overline{T_{j}}) \in \Tri \mk $ and
therefore $(T_{1},\overline{T_{2}},\overline{T_{3}})\in \Tri \mk $.

If $De_{i}\neq 0$, then from equation~(\ref{eqn:J1}) with $i=j$,
\begin{equation*}
2e_{i}\cdot De_{i}=2De_{i}
\end{equation*} 
so $De_{i}$ is an eigenvector of the multiplication $L_{e_{i}}$ with
eigenvalue $1$, i.e. $De_{i} \in P_{j}(\mk )+P_{k}(\mk )$ where $(i,j,k)$
are distinct. Write
\begin{equation*}
	De_{i} = P_{j}(x_{ij} )+P_{k}(x_{ik} );
\end{equation*}
then equation~(\ref{eqn:J1}) with $i\neq j$ gives
\begin{equation*}
	e_{i}\cdot P_{k}(x_{jk})+e_{i}\cdot
P_{i}(x_{ji})+P_{j}(x_{ij})\cdot e_{j}+P_{k}(x_{ik})\cdot e_{j} =0.
\end{equation*}
Thus
\begin{equation*}
	P_{k}(x_{jk}+x_{ik})=0.
\end{equation*}
It follows that the action of any derivation on the $e_{i}$ must be of
the form $F_{1}(x)+F_{2}(y)+F_{3}(z)$ where
\begin{align}
\label{eqn:Fone}
	F_{i}(x)e_{i} &=0 \notag \\
	F_{i}(x)e_{j}=-&F_{i}(x)e_{k}=P_{i}(x),
\end{align}
$(i,j,k)$ being a cyclic permutation of $(1,2,3)$. Hence $\Der H_{3}(\mk
)\subseteq \Tri \mk \oplus \mk ^{3}$.

To show that such derivations $F_{i}(x)$ exist and therefore the
inclusion just mentioned is an equality, consider the
operation of commutation with the matrix
\begin{align*}
	X&=\begin{pmatrix} 0 & -z & \bar{y} \\ \bar{z} & 0 & -x \\ -y &
\bar{x} & 0 \end{pmatrix} \\
	&= X_{1}(x)+X_{2}(y)+X_{3}(z)
\end{align*}
i.e. define $F_{i}(x)=C_{X_{i}(x)}$ where $C_{X}:H_{3}(\mk ) \to H_{3}(\mk
)$ is the commutator map
\begin{equation}
\label{eqn:com}
	C_{X}(H)=XH-HX.
\end{equation}
This satisfies equation~(\ref{eqn:Fone}) and also
\begin{align}
\label{eqn:FonP}
	F_{i}(x)P_{i}(y) &= -2 ( x,y ) (e_{j}-e_{k}) \notag \\
	F_{i}(x)P_{j}(y) &= -P_{k}(\bar{y}\, \bar{x})  \\
	F_{i}(x)P_{k}(y) &= P_{j}(\bar{x}\, \bar{y}). \notag
\end{align}
It is a derivation of $H_{3}(\mk )$ by
virtue of the matrix identity
\begin{equation}
\label{eqn:M1}
	[X,\{ H,K \}] = \{ [X,H],K\} + \{ H,[X,K]\}
\end{equation}
(in which square brackets denote commutators and round brackets denote
anticommutators), which we will prove separately in lemma~\ref{lemma:Matt}.

The Lie brackets of these derivations follow from another matrix
identity which is also proved in lemma~\ref{lemma:Matt},
\begin{equation}
\label{eqn:M2}
	[X,[Y,H]]-[Y,[X,H]]=[[X,Y],H]-E(X,Y)H
\end{equation}
where $E(X,Y)\in \so (\mk ^{\prime})$ is defined by
\begin{equation*}
	E(X,Y)z= \sum _{ij} [x_{ij},y_{ji},z],
\end{equation*}
$x_{ij},y_{ji}$ being the matrix elements of $X$ and $Y$. If
$X=X_{i}(x)$ and $Y=X_{j}(y)$ we have $D(X,Y)=0$ and
\begin{equation*}
	[X_{i}(x),X_{j}(y)]=X_{k}(\bar{y}\bar{x})
\end{equation*}
where $(i,j,k)$ is a cyclic permutation of $(1,2,3)$. This yields the
Lie bracket~(\ref{eqn:2}). If $X=X_{i}(x)$ and $Y=X_{j}(y)$, the matrix commutator
$Z=[X,Y]$ is diagonal with $z_{ii}=0$, $z_{jj}=y\bar{x}-x\bar{y}$ and
$z_{kk}=\bar{y}x-\bar{x}y$ ($i,j,k$ cyclic). Hence the action of the
commutator $[F_{i}(x),F_{i}(y)]=C_{z}+E(X,Y)$ on $H_{3}(\mk )$ is
\begin{align*}
	[F_{i}(x),F_{i}(y)]e_{m}&= 0 \quad (m=i,j,k) \\
	[F_{i}(x),F_{i}(y)]P_{i}(w)&=
P_{i}(z_{jj}w-wz_{kk}-2[x,y,w])=P_{i}(T_{1}w) \\
	&=4P_{i}(S_{xy}w) \quad \text{by equation~(\ref{eqn:Sdef}).} \\
	[F_{i}(x),F_{i}(y)]P_{j}(w)&=P_{j}(z_{kk}w-2[x,y,w]) \\
	&=P_{j}(\bar{y}(xw)-\bar{x}(yw)) \\
	[F_{i}(x),F_{i}(y)]P_{k}(w) &= P_{k}(-wz_{jj}-2[x,y,w]) \\
	&= P_{k}((wx)\bar{y}-(wy)\bar{x}).
\end{align*} 
Comparing with lemma~\ref{lemma:3}, we see that
\begin{align*}
	[F_{i}(x),F_{i}(y)]P_{i}(w) &=
P_{i}(T_{1}w)=P_{i}(T_{i}^{\prime}w) \\
	[F_{i}(x),F_{i}(y)]P_{j}(w) &=
P_{j}(\overline{T_{2}}w)=P_{j}(\overline{T_{j}^{\prime}}w) \\
	[F_{i}(x),F_{i}(y)]P_{k}(w) &=
P_{k}(\overline{T_{3}}w)=P_{k}(\overline{T_{k}^{\prime}}w) 
\end{align*}
where $(T_{1},T_{2},T_{3})=T_{xy}$, so that $T^{\prime}=\theta
^{1-i}(T_{xy})$. This establishes the Lie bracket~(\ref{eqn:3}).
\end{proof}
\end{theorem}

The matrix identities needed in the proof of
Theorem~\ref{theorem:banana} are contained in the following, in which we
include a third identity for the sake of completeness:
\begin{lemma}
\label{lemma:composition}
Let $\mk$ be a composition algebra, let $H$, $K$ and $L$ be hermitian
$3\times 3$ matrices with entries from $\mk$, and let $X,Y$ be traceless
antihermitian matrices over $\mk$. Then
\begin{subequations}
\begin{align}
	[X,\{ H,K\} ] &= \{ [X,H],K\} + \{ H,[X,K]\} \\
	[X,[Y,H]]-[Y,[X,H]] &=[[X,Y],H] +E(X,Y)H \\
	\{ H \{ K,L\} \} - \{ K\{ H,L \} \} &= [[H,K],L]+E(H,K)L
\end{align}
\end{subequations}
where $E(X,Y)\in \so (\mk )$ is defined for any $3\times 3$ matrices
$X,Y$ by
\begin{equation}
	E(X,Y)z = \sum _{ij} [x_{ij},y_{ji},z].
\end{equation}
In these matrix identities the square brackets denote commutators and
the chain brackets denote anticommutators.
\begin{proof}
We consider first part (a). The difference between the two sides can be written in terms
of matrix associators, where the $(i,j)$th element is
\begin{multline}
\label{eqn:paddington}
	\sum _{mn}(
[x_{im},h_{mn},k_{nj}]+[x_{im},k_{mn},h_{nj}]\\
+[k_{im},h_{mn},x_{nj}]-[h_{im},x_{mn},k_{nj}]-[k_{im},x_{mn},h_{nj}])
 .
\end{multline}
Suppose $i\neq j$ and let $k$ be the third index. Since the diagonal
elements of $H$ and $K$ are real, any associator containing them
vanishes. Hence the terms containing $x_{ij}$ or $x_{ji}$ are
\begin{multline*}
	\sum_{n}(
[x_{ij},h_{jn},k_{nj}]+[x_{ij},k_{jn},h_{nj}]) + \sum_{m} (
[h_{im},k_{mi},x_{ij}]+ [k_{im},h_{mi},x_{ij}]) \\
-[h_{ij},x_{ji},k_{ij}]+[k_{ij},x_{ji},h_{ij}]=0
\end{multline*}
by the alternative law, the hermiticity of $H$ and $K$, and the fact
that an associator changes sign when one of its elements is conjugated.
The terms containing $x_{ik}$ or $x_{ki}$ are
\begin{equation*}
[x_{ik},h_{ki},k_{ij}]+[x_{ik},k_{ki},h_{ij}]-[h_{ik},x_{ki},k_{ij}]- 
[k_{ik},x_{ki},h_{ij}]=0
\end{equation*}
using also $x_{ki}=-\bar{x}_{ik}$. Similarly, the terms containing
$x_{jk}$ or $x_{kj}$ vanish. Finally, the terms containing
$x_{ii},x_{jj}$ and $x_{kk}$ are
\begin{multline*}
[x_{ii},h_{ik},k_{kj}]+[x_{ii},k_{ik},h_{kj}]+[h_{ik},k_{kj},x_{jj}]+[k_{ik},h_{kj},x_{jj}]
 \\ -[h_{ik},x_{kk},k_{kj}]-[k_{ik},x_{kk},h_{kj}] =0
\end{multline*}
since $x_{ii}+x_{jj}+x_{kk}=0$.

Now consider the $(i,i)$th element. The last two terms of
equation~(\ref{eqn:paddington}) become
\begin{equation*}
-\sum_{mn} \left( [h_{im},x_{mn},k_{ni}]+[k_{in},x_{nm},h_{mi}]\right)
=0.
\end{equation*}
Let $j$ be one of the other two indices. The terms containing $x_{ij}$
or $x_{ji}$ are 
\begin{equation*}
[x_{ij},h_{jk},k_{ki}]+[x_{ij},k_{jk},h_{ki}]+[h_{ik},k_{kj},x_{ji}]+[k_{ik},h_{kj},x_{ji}]=0,
\end{equation*}
where $k$ is the third index. There are no terms containing $x_{jk}$ or
$x_{kj}$. The terms containing $x_{ii},x_{jj}$ or $x_{kk}$ are
\begin{multline*}
\sum_{n} \left( [x_{ii},h_{in},k_{ni}]+[x_{ii},k_{in},h_{ni}]\right) \\
+\sum_{m} \left( [h_{im},k_{mi},x_{ii}]+[k_{im},h_{mi},x_{ii}]\right)=0.
\end{multline*}
Thus in all cases the expression~(\ref{eqn:paddington}) vanishes, proving (a).
Parts (b) and (c) are proved by similar arguments, which the reader will
find more entertaining to write than to read.
\end{proof}
\end{lemma}

In any Jordan algebra, the commutator of two multiplication operators
$L_{x}$ and $L_{y}$ is a derivation (this fact is used in the
construction of the magic square Lie algebras $L_{3}(\mk _{1},\mk
_{2})$; see ~\cite{Tits62}). In the case of the Jordan algebra
$H_{3}(\mk )$, we can identify these derivations as follows:
\begin{lemma}
\label{lemma:lion}
In $H_{3}(\mk )$, where $\mk $ is any composition algebra,
\begin{align*}
	[L_{e_{i}},L_{e_{j}}] &= 0 \\
	[L_{e_{i}},L_{P_{j}(x)}] &= \epsilon _{ij}L_{P_{j}(x)}
\end{align*}
where $\epsilon _{ij}=0$ if $i=j$, otherwise $\epsilon _{ij}$ is the
sign of the permutation $(i,j,k)$ of $(1,2,3)$ where $k$ is the third index,
\begin{align*}
	[L_{P_{i}(x)},L_{P_{i}(y)}] &= -T_{x,y} \\
	[L_{P_{i}(x)},L_{P_{j}(y)}] &= F_{k}(\bar{y}\, \bar{x})
\end{align*}
where $e_{i},P_{i}(x) \in H_{3}(\mk )$ are defined by~(\ref{eqn:Jdef}) and
$F_{k}(x) \in \Der H_{3}(\mk )$ is given by~(\ref{eqn:Fone})
and~(\ref{eqn:FonP}).
\begin{proof}
Straightforward calculation from~(\ref{eqn:J1}-~\ref{eqn:J4}).
\end{proof}
\end{lemma}

\noindent The proof of Theorem~\ref{theorem:banana} suggests an alternative
description of $\Der H_{3}(\mk )$, which leads us to identify it as $\sa
(3,\mk )$:
\begin{theorem}
\label{theorem:apple}
For any composition algebra $\mk$,
\begin{equation}
	\Der H_{3}(\mk ) = \Der \mk \ds A_{3}^{\prime}(\mk)
\end{equation}
in which $\Der \mk$ is a Lie subalgebra, the Lie brackets between $\Der
\mk $ and $A_{3}^{\prime}(\mk)$  are given by the elementwise action of 
$\Der \mk$ on $3\times 3$ matrices, and
\begin{equation*}
	[X,Y]=(XY-YX)^{\prime}+\tfrac{1}{3}D(X,Y)
\end{equation*}
where $X,Y \in A_{3}^{\prime}(\mk)$,
\begin{equation*}
	(XY-YX)^{\prime}=XY-YX-\tfrac{1}{3}\tr (XY-YX)\I \quad \in
A_{3}^{\prime}(\mk)
\end{equation*}
and
\begin{equation*}
	D(X,Y)=\sum_{ij} D(x_{ij},y_{ji}) \quad \in \Der \mk
\end{equation*}
$x_{ij},y_{ji}$ being the matrix elements of $X$ and $Y$ and $D(x,y)$
being the derivation $D_{x,y}$ defined in equation~(\ref{eqn:Ddef}).
\begin{proof}
By Lemma~\ref{lemma:pig} and Theorem~\ref{theorem:banana}
\begin{equation}
\Der H_{3}(\mk ) = \Der \mk \ds 2\mk ^{\prime} \ds 3\mk .
\end{equation}
Identify $(a,b)+(x,y,z) \in 2\mk ^{\prime} \ds 3\mk$ with the traceless
antihermitian matrix
\begin{equation*}
	X= \begin{pmatrix} -a-b & -z & \bar{y} \\ \bar{z} & a & -x \\ -y
& \bar{x} & b \end{pmatrix} \quad \in A_{3}^{\prime}(\mk);
\end{equation*}
then the actions of $2\mk ^{\prime}$ and $3\mk $ on $H_{3}(\mk )$
defined in Theorem~\ref{theorem:banana} are together equivalent to the
commutator action $C_{x}$ defined by equation~(\ref{eqn:com}). By
Lemma~\ref{lemma:composition}(b),
\begin{equation*}
	[C_{X},C_{Y}]= C_{(XY-YX)^{\prime}} + C_{t\I }+E(X,Y)
\end{equation*}
where
\begin{align*}
	t&= \tfrac{1}{3}\tr (XY-YX) \\
	&=\tfrac{1}{3}\sum_{ij}(x_{ij}y_{ji}-y_{ji}x_{ij}).
\end{align*}
Now $C_{t\I }+E(X,Y)$ acts elementwise on matrices in $H_{3}(\mk)$
according to the map $D:\mk \to \mk $ given by
\begin{align*}
	Dz&=[t,z]+E(X,Y)z \\
	&= \sum_{ij}\left(
\tfrac{1}{3}[[x_{ij},y_{ji}],z]-[x_{ij},y_{ji},z]\right) \\
	&=\tfrac{1}{3}D(X,Y)Z.
\end{align*}
Hence the bracket $[X,Y]$ is as stated.
\end{proof}
\end{theorem}

Finally we use Theorem~\ref{theorem:banana} to give a description of
$L_{3}(\mk _{1},\mk _{2})$ which makes manifest the symmetry between
$\mk_{1}$ and $\mk_{2}$. 
\begin{theorem}
\label{theorem:peach}
For any two composition algebras $\mk _{1},\mk _{2}$,
\begin{equation}
L_{3}(\mk _{1}, \mk _{2}) = \Tri \mk _{1}\oplus \Tri \mk _{2} \ds 3\mk
_{1} \otimes \mk _{2}
\end{equation}
in which $\Tri \mk _{1}\oplus \Tri \mk _{2}$ is a Lie subalgebra;
\begin{align}
	[T_{1},F_{i}(x\otimes y)] &= F_{i}(T_{1i}x_{1}\otimes x_{2})
\quad \in 3\mk_{1}\otimes \mk_{2} \label{eqn:th3.1} \\
	[T_{2},F_{i}(x\otimes y)] &= F_{i}(x_{1}\otimes T_{2i}x_{2}) \quad \in 3\mk_{1}\otimes
\mk_{2} \label{eqn:th3.2} 
\end{align}
if $T_{\alpha }=(T_{\alpha 1},\overline{T}_{\alpha 2}, \overline{T}_{\alpha
3})\, \in \Tri \mk \, (\alpha = 1,2)$, and 
\begin{align}
F_{1}(x_{1}\otimes x_{2}) + F_{2}(y_{1}\otimes
y_{2}) + F_{3}(z_{1}\otimes z_{2})&=(x_{1}\otimes x_{2},y_{1}\otimes
y_{2},z_{1}\otimes z_{2}) \notag \\
&\hspace{2cm} \in 3\mk _{1}\otimes \mk _{2}; \notag\\
[F_{i}(x_{1}\otimes x_{2}) , F_{j}(y_{1}\otimes y_{2})] &=
F_{k}(\bar{y}_{1}\bar{x}_{1}\otimes \bar{y}_{2}\bar{x}_{2})
\label{eqn:th3.3} \\
&\hspace{2cm} \in 3\mk _{1} \otimes \mk _{2} \notag
\end{align}
if $x_{\alpha },y_{\alpha }\in \mk _{2}$ and $(i,j,k)$ is a cyclic
permutation of $(1,2,3)$; and
\begin{multline}
\label{eqn:th3.4}
[F_{i}(x_{1}\otimes x_{2}) , F_{i}(y_{1}\otimes y_{2})] = \langle
x_{2},y_{2}\rangle \theta ^{1-i}T_{x_{1}y_{1}}+\langle
x_{1},y_{1}\rangle \theta ^{1-i}T_{x_{2}y_{2}} \\
\hspace{4.5cm} \in \Tri \mk_{1} \oplus \Tri \mk_{2}
\end{multline}
\begin{proof}
We can write
\begin{equation*}
	H_{3}^{\prime}(\mk _{1})=2\mr \oplus 3\mk
\end{equation*}
by identifying $(\alpha,\beta )+(x,y,z)\in 2\mr \oplus 3\mk $ with the
matrix
\begin{align*}
&\begin{pmatrix} -\alpha -\beta & z & \bar{y} \\ \bar{z} & \alpha & x
\\ y & \bar{x} & \beta \end{pmatrix} \in H_{3}^{\prime}(\mk ), \\
&=\alpha (e_{2}-e_{1})+\beta (e_{3}-e_{1})+P_{1}(x)+P_{2}(y)+P_{3}(z)
\end{align*}
in the notation of theorem~\ref{theorem:banana}.
Then the vector space structure~(\ref{eqn:L3a}) of $L_{3}(\mk _{1},\mk_{2})$
can be written using Theorem~\ref{theorem:banana}, as
\begin{align*}
	L_{3}(\mk _{1},\mk_{2}) &= \Der H_{3}(\mk _{1}) \ds H_{3}^{\prime
}(\mk _{1})\otimes \mk _{2}^{\prime} \ds \Der \mk_{2} \\
	&= (\Tri \mk \ds 3\mk _{1})\ds (2\mk_{2}^{\prime} \ds
3\mk_{1}\otimes 2\mk _{2}^{\prime}) \ds \Der \mk_{2} \\
	&= \Tri \mk _{1} \ds (\Der \mk _{2} \ds 2\mk _{2}^{\prime}) \ds
(3\mk_{1}\otimes \mk_{2}^{\prime} \ds 3\mk _{1}) \\
	&\cong \Tri \mk _{1} \ds \Tri \mk _{2} \ds 3\mk _{1}\otimes \mk
_{2}.
\end{align*}

We use the following notation for the elements of the five subspaces of
$L_{3}(\mk _{1},\mk _{2})$:
\begin{enumerate}
\item $\Tri \mk \subset \Der H_{3}(\mk _{1})$ contains elements
$T=(T_{1},\overline{T}_{2},\overline{T}_{3})$ acting on
$H_{3}^{\prime}(\mk _{1})$ as in Theorem~\ref{theorem:banana}:
\begin{equation*}
Te_{i}=0, \quad TP_{i}(x)=P_{i}(T_{i}x) \quad (x\in \mk ; i = 1,2,3)
\end{equation*}
\item $3\mk _{1}$ is the subspace of $\Der H_{3}(\mk _{1})$ containing
the elements $F_{i}(x)$ defined in Theorem~\ref{theorem:banana}; these
will be identified with the elements $F_{i}(x\otimes 1)\in 3\mk
_{1}\otimes \mk _{2}^{\prime}$.
\item $2\mk _{2}^{\prime}$ is the subspace $\Delta \otimes \mk
_{2}^{\prime}$ of $H_{3}(\mk _{1})\otimes \mk _{2}^{\prime}$, where
$\Delta \subset H_{3}^{\prime}(\mk _{1})$ is the subspace of
real, diagonal, traceless matrices and is identified with the subspace of
$\Tri \mk $ as described in Lemma~\ref{lemma:pig}. We will regard $2\mk
_{2}^{\prime}$ as a subspace of $3\mk
_{2}^{\prime}$, namely
\begin{equation*}
	2\mk _{2}^{\prime}= \{ (a_{1},a_{2},a_{3})\in 3\mk _{2}^{\prime}
: a_{1}+a_{2}+a_{3}=0 \}
\end{equation*}
and identify $\mathbf{a}=(a_{1},a_{2},a_{3})$ with the $3\times 3$
matrix
\begin{equation*}
	\Delta (\mathbf{a})= \begin{pmatrix} a_{1} & 0 & 0 \\ 0 & a_{2}
& 0 \\ 0 & 0 & a_{3} \end{pmatrix} \in H_{3}^{\prime}(\mk _{1})\otimes
\mk _{2}^{\prime}
\end{equation*}
and with the triality
$T(\mathbf{a})=(T_{1},\overline{T}_{2},\overline{T}_{3})$ where
$T_{i}=L_{a_{j}}-R_{a_{k}}$ (see the remark after the proof of
Lemma~\ref{lemma:splodge}).
\item $3\mk _{1}\otimes \mk _{2}^{\prime}$ is the subspaces of
$H_{3}(\mk _{1}\otimes \mk _{2}^{\prime}$ spanned by elements
$P_{i}(x)\otimes a$ $(i=1,2,3 : x \in \mk _{1}, a \in \mk
_{2}^{\prime})$; it is also a subspace of $3\mk _{1}\otimes \mk _{2}$ in
the obvious way.
\item $\Der \mk _{2}$ is a subspace of $\Tri \mk _{2}$, a derivation $D$
being identified with $(D,D,D) \in \Tri \mk _{2}$.
\end{enumerate}

To complete the proof we must verify that the Lie brackets defined by
Tits (see Section~\ref{sec:ms3}) coincide with those in the statement of
the theorem. The above decomposition of $L_{3}(\mk_{1},\mk_{2})$ into
five parts gives us fifteen types of bracket to examine. We will write
$[\, ,\, ]_{\text{Tits}}$ for the bracket defined in section~\ref{sec:ms3} and
$[\, ,\, ]_{\text{here}}$ for that defined above.
\begin{enumerate} 
\item $[\Tri \mk_{1}, \Tri \mk _{2}]$: For $T_{1},T_{2} \in \Tri \mk
_{1}$, $[T_{1},T_{2}]_{\text{Tits}}$ is the bracket in $\Der H_{3}(\mk _{1})$,
which by theorem~\ref{theorem:banana} is the same as $[T_{1},T_{2}]_{\text{here}}$.
\item $[\Tri \mk _{1},3\mk _{1}]$: For $T \in \Tri \mk _{1}$,
$F_{i}(x\otimes 1)\in 3\mk _{1}$,
\begin{align*}
	[T,F_{i}(x)]_{\text{Tits}} &= F_{1}(T_{i}x) \ \text{see Theorem~\ref{theorem:banana}} \\
	&=F_{i}(T_{i}x \otimes 1) = [T,F_{i}(x \otimes 1)]_{\text{here}}.
\end{align*}
\item $[\Tri \mk _{1}, 2\mk _{2}]$: For $T_{1} \in \Tri \mk _{1}$,
$(a,b,c)\in 2\mk _{2}^{\prime}$,
\begin{equation*}
	[T_{1},(a,b,c)]_{\text{Tits}}=[T, e_{1}\otimes a + e_{2}\otimes b +
e_{3} \otimes c] = 0
\end{equation*}
since in Theorem~\ref{theorem:banana} $\Tri \mk _{1}$ was obtained as the
subspace of derivations which annihilate the diagonal matrices $e_{i}$.
On the other hand,
\begin{equation*}
	[T_{1}, (a,b,c)]_{\text{here}} = [T_{1},T_{2}(a,b,c)]=0.
\end{equation*}
\item $[\Tri \mk _{1}, 3\mk _{1}\otimes \mk _{2}^{\prime}]$: For $T_{1}
\in \Tri \mk _{1}$, $P_{i}(x \otimes a) \in 3\mk _{1} \otimes \mk
_{2}^{\prime}$,
\begin{equation}
	[T_{1}, P_{i}(x\otimes a)=P_{i}(T_{1i}x\otimes a)=[T_{1}, P_{i}(x \otimes a)]_{\text{here}}.
\end{equation}
\item $[\Tri \mk _{1}, \Der \mk _{2}]_{\text{Tits}} \subset [\Der H_{3}(\mk
_{1}), \Der \mk _{2}]=0$, while \\ \mbox{$[\Tri \mk _{1}, \Der \mk _{2}]_{\text{here}}
\subset [\Tri \mk _{1}, \Tri \mk _{2}]=0$.}
\item $[3\mk _{1}, 3\mk _{1}]$: $3\mk _{1} = 3\mk _{1}\otimes \mr$ is
spanned by $ F_{i}(x) = F_{i}(x\otimes 1)$ \\ $(i=1,2,3; x \in \mk _{1})$,
and $[F_{i}(x),F_{j}(y)]_{\text{Tits}}$ is given by
Theorem~\ref{theorem:banana}, while $[F_{i}(x\otimes 1), F_{j}(y \otimes 1)]_{\text{here}}$ is the same since
$T_{x_{2},y_{2}}=0$ if $x_{2},y_{2} \in \mr$.
\item $[3\mk _{1}, 2\mk _{2}^{\prime}]$: For $F_{i}(x) \in 3\mk _{1}$, $
\mathbf{a}=(a_{1},a_{2},a_{3}) \in 2\mk _{2}^{\prime}$ with 
$(a_{1}+a_{2}+a_{3}=0)$,
\begin{align*}
	[F_{i}(x),\mathbf{a}]_{\text{Tits}} &= [F_{i}(x), \sum e_{i}\otimes
a_{i}] \in [\Der H_{3}(\mk _{1}), H_{3}^{\prime}(\mk _{1}) \otimes \mk
_{2}^{\prime}] \\
	&= P_{i}(x) \otimes (a_{j}-a_{k}) \hspace{4cm} \text{by~(\ref{eqn:Fone})}
\end{align*}
while
\begin{align*}
	[F_{i}(x),\mathbf{a}]_{\text{here}} &= [F_{i}(x\otimes 1),
T(\mathbf{a})] \ \in [3\mk_{1}\otimes \mk _{2}, \Tri \mk _{2}] \\
	&= F_{i}(x \otimes (a_{j}-a_{k})).
\end{align*}
\item $[3\mk_{1}, 3\mk_{1} \otimes \mk _{2}^{\prime}]$: For $F_{i}(x) \in
3\mk _{1}, F_{j}(y\otimes a) \in 3\mk_{1} \otimes \mk _{2}^{\prime}$,
\begin{align*}
	[F_{i}(x), F_{i}(y\otimes a)]_{\text{Tits}} &= [F_{i}(x), P_{i}(y)
\otimes a ] \ \in [\Der H_{3}(\mk_{1}), H_{3}^{\prime}(\mk _{1}) \otimes
\mk _{2}^{\prime}] \\
	&= -2 \langle x,y\rangle (e_{j}-e_{k})\otimes a \ \in
H_{3}^{\prime}(\mk _{1}\otimes \mk _{2}^{\prime} \\
	&= -2 \langle x,y \rangle (a_{1},a_{2}, a_{3}) \ \in 2\mk
_{2}^{\prime}
\end{align*}
where $a_{i}=0, a_{j}=a, a_{k}=-a$ ($i,j,k$ cyclic). On the other hand
\begin{align*}
	[F_{i}(x), F_{i}(y\otimes a)]_{\text{here}} &= [F_{i}(x\otimes 1),
F_{i}(y\otimes a)] \\
	&=-\langle x,y \rangle \theta ^{1-i}T_{1,a} \\
	&=[F_{i}(x),F_{i}(y\otimes a)]_{\text{Tits}}
\end{align*}
since $T_{1,a}=(2L_{a}+R_{a},2R_{a},2L_{a})$ which is identified with
$(0,2a,-2a)$ \\$ \in \mk _{2}^{\prime}$ in paragraph 3 above. If $i \neq j$
and $(i,j,k)$ is a cyclic permutation of $1,2,3)$,
\begin{align*}
	[F_{i}(x),F_{j}(y\otimes a)]_{\text{Tits}} &= [F_{i}(x),P_{j}(y\otimes
a)] \ \in [\Der H_{3}(\mk _{1}), H_{3}^{\prime}(\mk _{1})\otimes \mk_{2}]
\\
	&= -P_{k}(\bar{y}\bar{x})\otimes a \in H_{3}^{\prime}(\mk _{1})
\otimes \mk _{2}^{\prime} \\
	&= -F_{k}(\bar{y}\bar{x})\otimes a \in 3\mk _{1} \otimes
\mk_{2}^{\prime}
\end{align*}
while $[F_{i}(x),F_{j}(y\otimes a)]_{\text{here}} = -F_{k}(\bar{y}\bar{x})
\otimes a$ since $\bar{a}=-a$. Similarly,
\begin{equation*}
	[F_{i}(x),F_{k}(y\otimes a)]_{\text{Tits}}=
F_{k}(\bar{x}\bar{y})\otimes a = [F_{i}(x),F_{k}(y\otimes a)]_{\text{here}}.
\end{equation*}
\item $[3\mk _{1}, \Der \mk _{2}]_{\text{Tits}} \in [\Der H_{3}(\mk _{1}), \Der
\mk _{2}]=0$ \\ and 
\begin{equation*}
[3\mk _{1},\Der \mk _{2}]_{\text{here}} \in [3\mk _{1}\otimes
\mr , \Der \mk _{2}]=0.
\end{equation*}
\item $[2\mk _{2}^{\prime},2\mk _{2}^{\prime}]$: For
$\mathbf{a},\mathbf{b} \in 2\mk _{2}^{\prime}$, with
$\mathbf{a}=(a_{1},a_{2},a_{3})$ and $\mathbf{b}=(b_{1},b_{2},b_{3})$
where $a_{1}+a_{2}+a_{3}=b_{1}+b_{2}+b_{3}=0$,
\begin{align*}
	[\mathbf{a},\mathbf{b}]_{\text{Tits}} &= [\sum e_{i} \otimes a_{i},
\sum e_{j}\otimes b_{j}] \\
	&= \sum_{i,j} \left( \langle e_{i},e_{j}\rangle
D_{a_{i},b_{j}}+(e_{i} \ast e_{j})\otimes \Im (a_{i}b_{j})+\langle
a_{i},b_{j}\rangle [L_{e_{i}},L_{e_{j}}]\right) \\
	&=\sum _{ij} \left( 2\delta
_{ij}D_{a_{i},b_{j}}+\delta_{ij}(2e_{i}-\tfrac{2}{3}\I ) \otimes
\tfrac{1}{2}[a_{i},b_{j}] \right) \\
	&=\sum_{i} \left( 2D_{a_{i},b_{i}}+\tfrac{1}{3}e_{i}\otimes
(2[a_{i},b_{i}]-[a_{j},b_{j}]-[a_{k},b_{k}])\right) \\
	&= [\mathbf{a},\mathbf{b}]_{\text{here}} \qquad
\end{align*}
\item $[2\mk _{2}^{\prime},3\mk _{1}\otimes \mk _{2}^{\prime}]$: For
$\mathbf{a}=(a_{1},a_{2},a_{3}) \in 2\mk _{2}^{\prime}$ and $
F_{i}(x\otimes b) \in 3\mk_{1}\otimes \mk _{2}^{\prime}$,
\begin{align*}
	[\mathbf{a}, F_{i}(x\otimes b)]_{\text{Tits}} &= [e_{i}\otimes
a_{i}+e_{j}\otimes a_{j}+e_{k}\otimes a_{k}, P_{i}(x)\otimes b] \\
	&= P_{i}(x)\otimes \tfrac{1}{2}[a_{j},b]+P_{i}(x)\otimes
\tfrac{1}{2}[a_{k},b]- \langle a_{j}-a_{k},b\rangle F_{i}(x)
\end{align*}
using Lemma~\ref{lemma:lion}. The first two terms belong to the subspace
$3\mk _{1}\otimes \mk _{2}^{\prime}$ of $H_{3}^{\prime}\otimes \mk
_{2}^{\prime}$ and the third to the subspace $3\mk _{1}$ of $\Der
H_{3}(\mk _{1})$, so together they constitute an element of $3\mk
_{1}\otimes \mk _{2}$:
\begin{align*}
	[\mathbf{a},F_{i}(x\otimes b)]_{\text{Tits}} &= F_{i}(x\otimes
(\tfrac{1}{2}[a_{j},b] = \tfrac{1}{2}[a_{k},b] - \langle a_{j}-a_{k},b
\rangle )) \\
	&= F_{i}(x\otimes (a_{j}b-ba_{k})) \\
	&= F_{i}(x\otimes T(\mathbf{a})_{i}b) \\
	&= [\mathbf{a},F_{i}(x\otimes b)]_{\text{here}}
\end{align*}
\item $[2\mk _{2}^{\prime},\Der \mk _{2}]$: Tits's bracket~(\ref{eqn:Tits})
coincides with the bracket in $\Tri \mk _{2}$ as given by
Lemma~\ref{lemma:pig}.
\item $[3\mk _{1}\otimes \mk _{2}^{\prime}, 3\mk _{1}\otimes \mk
_{2}^{\prime}]$: For $P_{i}(x), P_{i}(y) \in 3\mk _{1}$ and $a,b \in \mk
_{2}^{\prime}$, if $i\neq j$ and $(i,j,k)$ is a cyclic permutation of
$(1,2,3)$ then
\begin{equation*}
	[P_{i}(x)\otimes a,P_{j}(y)\otimes
b]_{\text{Tits}}=P_{k}(\bar{y}\bar{x})\otimes \tfrac{1}{2}[a,b]-\langle
a,b \rangle F_{k}(\bar{y}\bar{x})
\end{equation*}
by equation~(\ref{eqn:J4}) and Lemma~\ref{lemma:lion}. This is an
element of $3\mk _{1}\otimes \mk _{2}^{\prime} \ds 3\mk _{1}\otimes \mr
$ which is identified with the following element of $3\mk _{1} \otimes
\mk _{2}$:
\begin{align*}
	F_{k}(\bar{y}\bar{x} \otimes \tfrac{1}{2} \left(
[a,b]-(a\bar{b}+b\bar{a})\right) &= F_{k}(\bar{y}\bar{x} \otimes
\bar{b}\bar{a}) \\
	&=[F_{i}(x\otimes a),F_{j}(y\otimes b)]_{\text{here}}.
\end{align*}
If $i=j$, then
\begin{multline*}
	[P_{i}(x)\otimes a, P_{i}(y)\otimes b]_{\text{Tits}} = 4\langle x,y
\rangle D_{a,b} +\\ \tfrac{1}{2}\langle x,y \rangle
(-2e_{i}+e_{j}+e_{k})\otimes [a,b]  -\langle a,b \rangle \theta
^{1-i}T_{x,y}
\end{multline*}
by Lemma~\ref{lemma:lion}. The second term belongs to the subspace $2\mk
_{2}^{\prime}$ and it is to be identified with the triality
$\tfrac{1}{3}\langle x,y \rangle \theta ^{1-i}T$ where
$T_{1}=\tfrac{1}{3}(L_{[a,b]}-R_{[a,b]})$,
$\overline{T}_{2}=-\tfrac{1}{3}(R_{[a,b]}+2L_{[a,b]})$ and
$\overline{T}_{3}=\tfrac{1}{3}(2R_{[a,b]}+L_{[a,b]})$. By
~(\ref{eqn:Teq2}), $T=T_{a,b}$. Hence 
\begin{equation*}
	[P_{i}(x)\otimes a, P_{i}(y)\otimes b]_{\text{Tits}}=\langle
x,y\rangle \theta ^{1-i}T_{a,b}-\langle a,b \rangle \theta
^{1-i}T_{x,y}.
\end{equation*}
\item $[3\mk _{1}\otimes \mk _{2}^{\prime}, \Der \mk _{2}]$ is given by
the action of $\Der \mk _{2}$ on the second factor of the tensor product
in both cases.
\item $[\Der \mk _{2}, \Der \mk _{2}]$ is given by the Lie bracket of
$\Der \mk _{2}$ in both cases.
\end{enumerate}
\end{proof}
\end{theorem}

\section{Magic Squares of $2\times 2$ Matrix Algebras: Proofs}
\label{sec:proofs2}
In this section we prove the following theorems.
\begin{theorem}
\label{theorem:cabbage}
For $\mk = \mr ,\mc ,\mh $ and $\mo $, 
\begin{align*}
	L_{2}(\mk _{1},\mk _{2}) &\cong \so (\mk _{1} \ds \mk _{2})\\
	L_{2}(\mk _{1},\mks _{2}) &\cong \so(\tfrac{1}{2}(\nu _{1} 
+\nu _{2}),\tfrac{1}{2}\nu _{2}).
\end{align*}
\end{theorem}
\begin{theorem}
\label{theorem:parsnip}
The following isomorphisms are true for $\mk = \mr ,\mc ,\mh $ and $\mo
$.
\begin{subequations}
\begin{align}
	L_{2}(\mk ,\mr )&\cong \Der H_{2}(\mk ) \label{eqn:parsnip1}\\
	L_{2}(\mk ,\mcs )&\cong \Str H_{2}(\mk ) \label{eqn:parsnip2}\\
	L_{2}(\mk ,\mhs )&\cong \Con H_{2}(\mk ) \label{eqn:parsnip3}
\end{align}
\end{subequations}
\end{theorem}

We prove Theorem~\ref{theorem:cabbage} by first showing that $L_2(\mk_{1} ,\mk_{2})$ is isomorphic to the Lie
algebra of the pseudo-orthogonal group $\O (\mk_{1}\ds\mk_{2})$ of
linear transformations of $\mk_{1}\ds \mk_{2}$ preserving the
quadratic form 
\begin{equation*}
	\left| x_{1}+x_{2}\right| ^{2} = \left| x_{1} \right| ^{2} + 
 \left| x_{2} \right| ^{2}. \quad (x_{1} \in \mk _{1}, x_{2} \in \mk _{2})
\end{equation*}
We then prove equation~(\ref{eqn:parsnip3})
and notice that the proof of this contains the isomorphisms for
equations~(\ref{eqn:parsnip1}) and~(\ref{eqn:parsnip2}).
The proofs of these equations will require the use of the following
Theorem, and associated Lemmas.
\begin{theorem}
\label{theorem:Nick}
The derivation algebra of $H_{2}(\mk )$ can be expressed in the form
\begin{equation}
\label{eqn:der}
	\Der H_{2}(\mk)=A_{2}^{\prime}(\mk)\ds \so(\mk^{\prime}).
\end{equation}
\end{theorem}
Our proof requires the use of the lemma
\begin{lemma}
\label{lemma:Matt}
Let $A\in A_{n}(\mk)$ and $X,Y \in H_{n}(\mk)$ where $\mk$ is any
alternative algebra. The identity
\begin{equation}
\label{eqn:ada}
	[A,\{ X,Y\} ]=\{ [A,X],Y\} +\{ X,[A,Y]\} ,
\end{equation}
(where the brackets denote commutators of matrices) holds if $n=2$ or if
$n=3$ and $\tr A=0$.
\begin{proof}
The $3\times 3$ case was proved in Lemma~\ref{lemma:composition}. The
$2\times 2$ case can be deduced from it by considering the $3\times 3$
matrices
\begin{equation*}
	\widetilde{A}=\begin{pmatrix} A & 0 \\ 0 & -\tr A \end{pmatrix},
\quad \widetilde{X}=\begin{pmatrix} X & 0 \\ 0 & 1 \end{pmatrix}, \quad
\widetilde{Y}=\begin{pmatrix} Y & 0 \\ 0 & 1 \end{pmatrix}.
\end{equation*}
\end{proof}
\end{lemma}

\begin{proof}[Proof of theorem~\ref{theorem:Nick}]
From lemma~\ref{lemma:Matt} we see that for each $A\in A_{2}(\mk)$ there
is a derivation $D(A)$ of $H_{2}(\mk)$ given by
\begin{equation*}
	D(A)(X) = AX-XA.
\end{equation*}

We consider $H_{2}(\mk)$ as a Jordan algebra with product
\begin{equation*}
	\begin{pmatrix} \alpha & x \\ \bar{x} & \beta \end{pmatrix} \cdot 
\begin{pmatrix} \gamma & y \\ \bar{y} & \delta
\end{pmatrix}=\begin{pmatrix} 2\alpha \gamma +2\Re (x\bar{y}) & (\gamma
+\delta )x+(\alpha +\beta )y \\ (\gamma
+\delta )\bar{x}+(\alpha +\beta )\bar{y} & 2\beta \delta +2\Re
(\bar{x}y) \end{pmatrix}.
\end{equation*}
We can write a matrix $A\in H_{2}(\mk)$ as follows
\begin{equation*}
	\begin{pmatrix} \alpha & x \\ \bar{x} & \beta \end{pmatrix} =
\lambda I+\mu E+P(x)
\end{equation*}
where $\lambda =\tfrac{1}{2}(\alpha +\beta )$, $\mu =\tfrac{1}{2}(\alpha -\beta
)$, $E=\begin{pmatrix} 1 & 0 \\ 0 & -1 \end{pmatrix}$ and
$P(x)=\begin{pmatrix} 0 & x \\ \bar{x} & 0 \end{pmatrix}$. Then the
Jordan multiplication can be rewritten as 
\begin{align*}
	E\cdot E &=I \\
	P(x) \cdot P(y) &= 2\langle x,y \rangle I \\
	E\cdot P(x) &= 0.
\end{align*}
Thus $H_{2}(\mk)$ can be identified with $\mj (V)$, the Jordan algebra
associated with the inner product space  $V=\mk \oplus \mr $. $\mj (V)$
is a subalgebra of
the anticommutator algebra of $\Cl (V)$, where $\mathbf{v} \cdot \mathbf{w}
= \langle \mathbf{v},\mathbf{w} \rangle 1$.
Derivations of this algebra must satisfy
\begin{align*}
	D(1)&=0 \\
	\langle 1, D(\mathbf{v})\rangle &=0.
\end{align*}
Thus 
\begin{equation*}
	\langle D(\mathbf{v}) , \mathbf{w} \rangle + \langle
D(\mathbf{w}), \mathbf{v} \rangle =0
\end{equation*}
i.e. $D$ is an antisymmetric map of $\mathbf{v}$. Hence $\Der
H_{2}(\mk) = \o (\mk \ds \mr)$.

Considering the matrix structure of $\o (\mk \ds \mr)$ we can write
this as $\o (\mk ) \ds \mk $. Consider the action of $\mk$ on the
$(\nu +1)\times 1$ column vectors $\begin{pmatrix} 0 \\1 \end{pmatrix}$
and $\begin{pmatrix} x \\0 \end{pmatrix}$. We express $k\in \mk $ as the final row
and column in a $(\nu +1)\times (\nu +1)$ block matrix. Then
\begin{align*}
	\begin{pmatrix} 0 & k \\ -k^{t} & 0 \end{pmatrix}
\begin{pmatrix} 0 \\ 1 \end{pmatrix} &= \begin{pmatrix} k \\ 0
\end{pmatrix}  \\
\begin{pmatrix} 0 & k \\ -k^{t} & 0 \end{pmatrix}
\begin{pmatrix} x \\ 0 \end{pmatrix} &= \begin{pmatrix} 0 \\ -k^{t}x 
\end{pmatrix} 
\end{align*}
Thus  $k$ maps $E$ to $P(k)$ and $P(x)$ to $-\langle k,x
\rangle E$. Now 
\begin{equation*}
	\left[ \begin{pmatrix} 0 & k \\ -\bar{k} & 0 \end{pmatrix},\begin{pmatrix} 0 & x \\ \bar{x} & 0
\end{pmatrix} \right] = 2\begin{pmatrix} \langle k,x \rangle  & -k \\
\bar{k} & -\langle k,x \rangle \end{pmatrix}
\end{equation*}
i.e. multiplication by $\begin{pmatrix} 0 & k \\ -k^{t} & 0
\end{pmatrix}$ in $\mj (V)$ is equivalent to  commutation with
 $\begin{pmatrix} 0 & -\tfrac{k}{2} \\ \tfrac{\bar{k}}{2} & 0
\end{pmatrix}$ in $H_{2}(\mk)$. 

We can split $\o (\mk)=\o (\mk ^{\prime} \oplus \mr)$ into $\o (\mk ^{\prime}) \ds \mk$.
Consider the
action of the $\nu \times \nu$ matrix $\begin{pmatrix} 0 & l \\ -l^{t} & 0
\end{pmatrix}$ with $l \in \mk^{\prime}$ on the vectors $\begin{pmatrix} y \\ 0
\end{pmatrix}$ and $\begin{pmatrix} 0 \\ 1 \end{pmatrix}$ with $y \in
\mk^{\prime}$:
\begin{align*}
	 \begin{pmatrix} 0 & l \\ -l^{t} & 0 \end{pmatrix}\begin{pmatrix} y \\ 0
\end{pmatrix} &= \begin{pmatrix} -l^{t}x \\ 0 \end{pmatrix} \\
	\begin{pmatrix} 0 & l \\ -l^{t} & 0 \end{pmatrix}\begin{pmatrix} 0 \\ 1
\end{pmatrix} &= \begin{pmatrix} 0 \\ l \end{pmatrix} 
\end{align*}
and we obtain (by a similar method) that multiplication by $\begin{pmatrix} 0 & l \\ -l^{t} & 0
\end{pmatrix}$ in $\mj (V)$ is equivalent to commutation with
$\begin{pmatrix} \tfrac{l}{2} & 0 \\ 0 & -\tfrac{l}{2} \end{pmatrix}$ in
$H_{2}(\mk )$. 
Further $\so(\mk ^{\prime})$ acts in $\mj (V)$ precisely as it does in
$H_{2}(\mk)$. Thus we have 
\begin{equation*}
	\Der H_{2}(\mk) = A_{2}^{\prime}(\mk ) \ds \so(\mk^{\prime})
\end{equation*}
as required. 
\end{proof}
The brackets in $A_{2}^{\prime}(\mk _{1})\ds \so (\mk_{1})$ are given by
\begin{align*}
	[A,A^{\prime}]&=AA^{\prime}-A^{\prime}A \\
	[S,A]&=S(A) \\
	[S,S^{\prime}]&=SS^{\prime}-S^{\prime}S
\end{align*}
with $A,A^{\prime} \in A_{2}^{\prime}(\mk_{1})$ and $S,S^{\prime}\in
\so(\mk_{1})$ and $S(A)$ describes $S$ acting elementwise on
$A$. When they arise in calculations we consider multiples of the
$2\times 2$ identity matrix $I_{2}$ to be elements of $\so(\mk_{1})$.
Finally, the following lemma holds.
\begin{lemma}
\label{lemma:Mum}
The Jacobi identity 
\begin{equation*}
	[ A,[B,H]]+[B,[H,A]]+[H,[A,B]]=0
\end{equation*}
holds for $A,B\in A_{2}^{\prime}(\mk)$ and $X\in H_{2}(\mk)$.
\end{lemma}
\begin{proof}[Proof of Theorem~\ref{theorem:cabbage}]
Using Theorem~\ref{theorem:Nick} we can write $L_{2}(\mk _{1}, \mk _{2})$ as
\begin{equation*}
	L_{2}(\mk _{1},\mk _{2}) =A_{2}^{\prime}(\mk _{1})\ds \so (\mk_{1}) \ds H^{\prime }_{2}(\mk
_{1})\otimes \mk ^{\prime }_{2}\ds \so (\mk ^{\prime }_{2}).
\end{equation*}
This can be considered to contain the following elements
\begin{align}
\label{eqn:Ann}
	J&=\begin{pmatrix}
		0 & 1 \\
		-1 & 0 \\
	\end{pmatrix} \ \in A_{2} ^{\prime}(\mr) \notag \\ 
	A_{1}&=\begin{pmatrix}
		a_1 & 0 \\
		0 & -a_1 \\
	\end{pmatrix} \ \in A_{2} ^{\prime}(\mk _{1}^{\prime}) \notag \\
	A_{2}&=\begin{pmatrix}
		0 & a_{2} \\
		a_{2} & 0 \\
	\end{pmatrix} \ \in A_{2} ^{\prime}(\mk _{1}^{\prime}) \notag \\
	F&\in \so(\mk _{1}^{\prime}) \\
	S&\in \so(\mk _{2}^{\prime}) \notag \\
	B_{1}&=\begin{pmatrix}
		1 & 0 \\
		0 & -1 \\
	\end{pmatrix} \otimes b_{1}  \ \in H_{2} ^{\prime}(\mr
)\otimes \mk _{2}^{\prime} \notag \\ 
	B_{2}&=\begin{pmatrix}
		0 & 1 \\
		1 & 0 \\
	\end{pmatrix} \otimes b_{1} \ \in H_{2} ^{\prime}(\mr
)\otimes \mk _{2}^{\prime} \notag \\ 
	C&=\begin{pmatrix}
		0 & c_1 \\
		-c_1 & 0 \\
	\end{pmatrix} \otimes c_{2} \ \in H_{2} ^{\prime}(\mk
_{1}^{\prime})\otimes \mk _{2}^{\prime}  \notag
\end{align}
with $a_{1},b_{1},c_{1}\in \mk _{1}^{\prime}$ and $a_{2},b_{2},c_{2}\in \mk _{2}^{\prime}$.

We define $\varphi :L_{2}(\mk _1,\mk _2)\rightarrow \so (\nu _{1}+\nu
_{2})$ by
\begin{multline}
	\varphi (J+A_{1}+A_{2}+F+S+B_{1}+B_{2}+C)= \\
		\begin{pmatrix}
		S & 2Gb_{1} & 2Gb_{2} & 2Gc_{2}c_{1}^{t} \\
		-2b_{1}^{t}G & 0 & 2 & 2a_{2}^{t} \\
		-2b_{2}^{t}G & -2 & 0 & -2a_{1}^{t} \\
		-2c_{1}c_{2}^{t}G & -2a_{2} & 2a_{1} & F \\
		\end{pmatrix}
\end{multline}
where now elements of $\mk _{i} ^{\prime}$ are identified with column vectors
in $\mr ^{\nu _{i}-1}$ and $G$ is the metric matrix for $\mk _{2}^{\prime}$. In the
Euclidean case it is merely the identity matrix whereas in the non-Euclidean
case it is the diagonal matrix consisting of $(\tfrac{\nu_{2}}{2})$
positive $1$'s
and $(\tfrac{\nu_{2}}{2}-1)$ negative $1$'s. The order of the positive and
negative elements is determined by the choice of $\pm 1$ in
the Cayley-Dickson calculation for $\mk_{2}$.
We show that $\psi$ is a Lie algebra isomorphism by calculating the
multiplication tables for the Lie brackets between the elements listed
in $L_{2}(\mk _{1},\mk _{2})$. We then calculate the equivalent brackets
in $\so (\mk_{1}\oplus \mk _{2})$ and show that they are equivalent. The
relevant tables are found on the following two pages  and are obtained
simply by applying the stated Lie brackets for each algebra to these
basis elements.
\newpage
\begin{landscape}
\begin{center}
\tiny
\begin{tabular}{|p{0.7cm}||p{2.2cm}|p{2.2cm}|p{2.0cm}|p{2.0cm}|p{2.0cm}|p{2.2cm}|p{2.0cm}|p{2.0cm}|}
\hline
 & $J$ & $A_{1}$ & $A_{2}$ & $B_{1}$ & $B_{2}$ &
$C$ & $F$ & $S$ \\ \hline \hline
$J$ & $0$ & $-2\begin{pmatrix} 0 & a_{1} \\ a_{1} & 0 \end{pmatrix}$ &
$2\begin{pmatrix} a_{2} & 0 \\ 0 & -a_{2} \end{pmatrix}$ & $-2
\begin{pmatrix} 0 & 1 \\ 1 & 0 \end{pmatrix} \otimes b_{1}$ &
$2\begin{pmatrix} 1 & 0 \\ 0 & -1 \end{pmatrix} \otimes b_{2}$ & $0$ &
$0$ & $0$ \\ \hline
$A_{1}^{\prime}$ & $2\begin{pmatrix} 0 & a_{1}^{\prime} \\ a_{1}^{\prime} & 0
\end{pmatrix}$ & $2\Im (a_{1}a_{1}^{\prime})I$ & $-2\Re
(a_{1}^{\prime}a_{1})J$ & $0$ & $2\begin{pmatrix} 0 & a_{1}^{\prime} \\ a_{1}^{\prime} & 0
\end{pmatrix} \otimes b_{2}$ & $2\Re
(c_{1}a_{1}^{\prime})\times \begin{pmatrix} 0 & 1 \\ 1 & 0 \end{pmatrix}
\otimes c_{2}$ & $F(A_{1}^{\prime})$ & $0$ \\ \hline
$A_{2}^{\prime}$ & $-2\begin{pmatrix} a_{2}^{\prime} & 0 \\ 0 &
-a_{2}^{\prime} \end{pmatrix}$ & $ -2\Re (a_{1}a_{2}^{\prime})J$ & $2\Im
(a_{2}a_{2}^{\prime})I$ & $2\begin{pmatrix} 0 & -a_{2}^{\prime} \\ a_{2}^{\prime}
& 0 \end{pmatrix} \otimes b_{1}$ & $0$ & $ 2\Re
(a_{2}^{\prime}c_{1})\times \begin{pmatrix} 1 & 0 \\ 0 & -1
\end{pmatrix} \otimes c_{2}$ & $ F(A_{2}^{\prime}$ & $0$ \\ \hline
$B_{1}^{\prime}$ & $2\begin{pmatrix} 0 & 1 \\ 1 & 0 \end{pmatrix}
\otimes b_{1}^{\prime}$ & $0$ & $2\begin{pmatrix} 0 & a_{2} \\ a_{2} & 0
\end{pmatrix} \otimes b_{1}^{\prime}$ & $D_{b_{1}^{\prime},b_{1}}$ &
$-2(b_{1}^{\prime},b_{2})J$ &
$ -2(b_{1}^{\prime},c_{2})\times \begin{pmatrix} 0 & c_{1} \\ c_{1} & 0
\end{pmatrix}$ & $0$ & $ \begin{pmatrix} 1 & 0 \\ 0 & -1 \end{pmatrix}
\otimes S(b_{1}^{\prime})$ \\ \hline
$B_{2}^{\prime}$ & $ -2\begin{pmatrix} 1 & 0 \\ 0  & -1\end{pmatrix}
\otimes b_{2}^{\prime}$ & $ 2\begin{pmatrix} 0 -a_{1} \\ a_{1} & 0
\end{pmatrix} \otimes b_{2}^{\prime}$ & $0$ &
$2(b_{1},b_{2}^{\prime})J$ & $D_{b_{2}^{\prime},b_{2}}$ & $
-2(b_{2}^{\prime},c_{2})\times \begin{pmatrix} -c_{1} & 0 \\ 0 & c_{1}
\end{pmatrix}$ & $0$ & $\begin{pmatrix} 0 & 1 \\ 1 & 0 \end{pmatrix}
\otimes S(b_{2}^{\prime})$ \\ \hline
$C^{\prime}$ & $0$ & $-2\Re (c_{1}^{\prime}a_{1}) \times \begin{pmatrix} 0
& 1 \\ 1 & 0 \end{pmatrix} \otimes c_{2}^{\prime}$ & $ 2\Re
(a_{2}c_{1}^{\prime})\times \begin{pmatrix} 1 & 0 \\ 0 & -1
\end{pmatrix} \otimes c_{2}^{\prime}$ & $2(b_{1},c_{2}^{\prime}) \times
\begin{pmatrix} 0 & c_{1}^{\prime} \\ c_{1}^{\prime} & 0 \end{pmatrix}$
& $-2(b_{2},c_{2}^{\prime})\times \begin{pmatrix} c_{1}^{\prime} & 0 \\ 0
& -c_{1}^{\prime} \end{pmatrix}$ &
$(c_{1},c_{1}^{\prime})D_{c_{2}^{\prime},c_{2}}- 2
(c_{2},c_{2}^{\prime})\times \Im (c_{1}^{\prime}c_{1})I$ &
$F\begin{pmatrix} 0 & c_{1}^{\prime} \\ c_{1}^{\prime} & 0 \end{pmatrix}
\otimes c_{2}^{\prime}$ & $\begin{pmatrix} 0 & c_{1}^{\prime} \\
-c_{1}^{\prime} & 0 \end{pmatrix} \otimes S(c_{2}^{\prime})$ \\ \hline
$F^{\prime}$ & $0$ & $-F^{\prime}(A_{1})$ & $-F^{\prime}(A_{2})$ & $0$ &
$0$ & $ -F^{\prime}\begin{pmatrix} 0 & c_{1} \\ c_{1} & 0 \end{pmatrix}
\otimes c_{2}$ & $F^{\prime}F-FF^{\prime}$ & $0$ \\ \hline
$S^{\prime}$ & $0$ & $0$ & $0$ & $-\begin{pmatrix} 1 & 0 \\ 0 & -1
\end{pmatrix} \otimes S^{\prime}(b_{1})$ & $-\begin{pmatrix} 0 & 1 \\ 1
& 0 \end{pmatrix} \otimes S^{\prime}(b_{2})$ & $-\begin{pmatrix} 0 &
c_{1} \\ c_{1} & 0 \end{pmatrix} \otimes S^{\prime}(c_{2})$ & $0$ &
$S^{\prime}S-SS^{\prime}$ \\ \hline
\end{tabular}
\end{center}
\end{landscape}
\newpage
\begin{landscape}
\begin{center}
\tiny
\begin{tabular}{|p{0.8cm}||p{1.4cm}|p{2.0cm}|p{2.0cm}|p{2.2cm}|p{2.2cm}|p{2.5cm}|p{2.0cm}|p{2.0cm}|}
\hline
 & $\psi (J)$ & $\psi (A_{1})$ & $\psi (A_{2})$
& $\psi (B_{1})$ & $\psi (B_{2})$ & $\psi (C)$ & $\psi (F)$ & $\psi (S)$ \\
\hline \hline
$\psi (J^{\prime})$ & $0$ & $-4a_{1}^{t}G_{24}+4a_{1}G_{42}$ &
$-4a_{2}^{t}G_{34}+4a_{2}G_{43}$ & $-4b_{1}GG_{13}+4Gb_{1}^{t}G_{31}$ &
$4b_{2}GG_{12}-4Gb_{2}^{t}G_{21}$ & $0$ & $0$ & $0$ \\ \hline
$\psi (A_{1}^{\prime})$ & $ 4a_{1}^{\prime \,
t}G_{24}-4a_{1}^{\prime}G_{42}$ &
$(-4a_{1}^{\prime}a_{1}^{t}+4a_{1}a_{1}^{\prime \, t})G_{44}$ &
$-4a_{2}^{t}a_{1}^{\prime}G_{23}+4a_{1}^{\prime \, t}a_{2}G_{32}$ & $0$
& $4Gb_{2}a_{1}^{\prime \, t}G_{14}-4a_{1}^{\prime}b_{2}^{t}GG_{41}$ & $
-4Gc_{2}c_{1}^{t}a_{1}^{\prime}G_{13}+4a_{1}^{\prime \,
t}c_{1}c_{2}^{t}GG_{31}$ & $-2a_{1}^{\prime \,
t}FG_{34}+2F^{t}a_{1}^{\prime}G_{43}$ & $0$ \\ \hline
$\psi (A_{2}^{\prime})$ & $ 4a_{2}^{\prime \,
t}G_{34}-4a_{2}^{\prime}G_{43}$ & $ 4a_{2}^{\prime \,
t}a_{1}G_{23}-4a_{1}^{t}a_{2}^{\prime}G_{32}$ &
$(-4a_{2}^{\prime}a_{2}^{t}+4a_{2}a_{2}^{\prime \, t})G_{44}$ &
$4a_{2}^{\prime}b_{1}^{t}GG_{41}-4Gb_{1}a_{2}^{\prime \, t}G_{14}$ & $0$
& $4Gc_{2}c_{1}^{t}a_{2}^{\prime}G_{12}-4a_{2}^{\prime \,
t}c_{1}c_{2}^{t}GG_{21}$ & $2a_{2}^{\prime \, t}FG_{24}-2Fa_{2}G_{42}$ &
$0$ \\ \hline 
$\psi (B_{1}^{\prime})$ & $4b_{1}^{\prime}GG_{13}-4b_{1}^{\prime \,
t}GG_{31}$ & $0$ & $4Gb_{1}^{\prime}a_{2}^{t}G_{14}-4a_{2}b_{1}^{\prime
\, t}GG_{41}$ & $(-4b_{1}^{\prime}b_{1}^{t}+4b_{1}b_{1}^{\prime \,
t})G_{11}$ & $4b_{2}^{t}GGb_{1}^{\prime}G_{32}-4b_{1}^{\prime \,
t}GGb_{2}G_{23}$ & $
4Gc_{2}c_{1}^{t}Gb_{1}^{\prime}G_{42}-4b_{1}^{\prime \,
t}GGc_{2}c_{1}^{t}G_{24}$ & $0$ & $2b_{1}^{\prime \,
t}GS^{t}G_{21}-2SGb_{1}^{\prime}G_{12}$ \\ \hline
$\psi (B_{2}^{\prime})$ & $4b_{2}^{\prime \,
t}GG_{21}-4Gb_{2}^{\prime}G_{12}$ & $4a_{1}b_{2}^{\prime \,
t}GG_{41}-4b_{2}^{\prime}Ga_{1}^{t}G_{14}$ & $0$ &
$4b_{1}^{t}GGb_{2}^{\prime}G_{23}-4b_{2}^{\prime \, t}GGb_{1}G_{32}$ &
$4(Gb_{2}b_{2}^{\prime \, t}G-Gb_{2}^{\prime}b_{2}^{t}G)G_{11}$ &
$4c_{1}c_{2}^{t}GGb_{2}^{\prime}G_{43}-4b_{2}^{\prime \,
t}GGc_{2}c_{1}^{t}G_{34}$ & $0$ & $2b_{2}^{\prime \,
t}GS^{t}G_{31}-2Sb_{2}^{\prime}GG_{13}$ \\ \hline
$\psi (C^{\prime})$ & $0$ & $4Gc_{2}^{\prime}c_{1}^{\prime \,
t}a_{1}G_{13}-4a_{1}^{t}c_{1}^{\prime}c_{2}^{\prime \, t}GG_{31}$ &
$4a_{2}^{t}c_{1}^{\prime}c_{2}^{\prime \,
t}GG_{21}-4Gc_{2}^{\prime}c_{1}^{\prime \, t}a_{2}G_{12}$ &
$4b_{1}^{t}GGc_{2}^{\prime}c_{1}^{\prime \,
t}G_{24}-4Gc_{2}^{\prime}c_{1}^{\prime \, t}Gb_{1}G_{42}$ &
$4b_{2}^{t}GGc_{2}^{\prime}c_{1}^{\prime \,
t}G_{34}-4c_{1}^{\prime}c_{2}^{\prime\, t}GGb_{2}G_{43}$ &
$4(Gc_{2}c_{1}^{t}c_{1}^{\prime}c_{2}^{\prime \,
t}G-Gc_{2}^{\prime}c_{1}^{\prime \,
t}c_{1}c_{2}^{t}G)G_{11}+4(c_{1}c_{2}^{t}GGc_{2}^{\prime}c_{1}^{\prime
\, t}-c_{1}^{\prime}c_{2}^{\prime \, t}GGc_{2}c_{1}^{t})G_{44}$ &
$2Gc_{2}^{\prime}c_{1}^{\prime \,
t}FG_{14}-2Fc_{1}^{\prime}c_{2}^{\prime \, t}GG_{41}$ & $
2c_{1}^{\prime}c_{2}^{\prime \,
t}GS^{t}G_{41}-2SGc_{2}^{\prime}c_{1}^{\prime \, t}G_{14}$ \\ \hline
$\psi (F^{\prime})$ & $0$ & $2a_{1}^{t}F^{\prime}G_{34}-2F^{\prime \,
t}a_{1}G_{43}$ & $2F^{\prime \,
t}a_{2}G_{42}-2a_{2}^{t}F^{\prime}G_{24}$ & $0$ & $0$ &
$2F^{\prime}c_{1}c_{2}^{t}GG_{41}-2Gc_{2}c_{1}^{t}F^{\prime}G_{14}$ & $(
F^{\prime}F-FF^{\prime})G_{44}$ & $0$ \\ \hline
$\psi (S^{\prime})$ & $0$ & $0$ & $0$ &
$2S^{\prime}Gb_{1}G_{12}-2b_{1}^{t}GS^{\prime \, t}G_{21}$ &
$2S^{\prime}Gb_{2}G_{13}-2b_{2}^{t}GS^{\prime \, t}G_{31}$ &
$2S^{\prime}Gc_{2}c_{1}^{t}G_{14}-2c_{1}c_{2}^{t}GS^{\prime \, t}G_{41}$
& $0$ & $(S^{\prime}S-SS^{\prime})G_{11}$ \\ \hline
\end{tabular}
\end{center}
\end{landscape}

\newpage

Since Lie brackets are antisymmetric the tables show that there are 25
non-zero brackets to compare. We will use the following lemma.
\begin{lemma}
If $\mathbf{a}_{1}$, $\mathbf{a}_{2} \in \mr ^{\nu -1}$ are the vector
representations of the hypercomplex numbers $a_{1},a_{2} \in \mk
^{\prime}$ then
\begin{enumerate}
\item
$4(\mathbf{a}_{1}\mathbf{a}_{2}^{t}+\mathbf{a}_{2}\mathbf{a}_{1}^{t})
\in \so (\mk _{1} \ds \mk _{2})$ is equivalent to $2\Im
(a_{2}a_{1})I \in L_{2}(\mk _{1},\mk _{2})$,
\item $\mathbf{a}_{2}^{t}\mathbf{a}_{1}\in \so (\mk _{1} \ds \mk _{2})$
is equivalent to $(a_{1},a_{2})=-\Re (a_{1}a_{2}) \in L_{2}(\mk _{1},\mk _{2})$.
\end{enumerate}
\begin{proof}
We shall use the notation that if $a=a_{i}e_{i}$ is a hypercomplex number in $\mk
^{\prime}$ then $\mathbf{a}=(a_{i})$ is its vector representation as a
column vector in $\mr ^{\nu
-1}$ where $\nu$ is the dimension of $\mk $.
\begin{enumerate}
\item The $(i,j)$th component of the matrix
$(\mathbf{a}\mathbf{b}^{t}+\mathbf{b}\mathbf{a}^{t})$ is
the element $a_{i}b_{j}+b_{i}a_{j}$. If this multiplies a third
hypercomplex number, say $\mathbf{c}=(c_{k})$ then the $n$th element of
the resulting vector will be
\begin{equation*}
	\sum_{k=1}^{\nu -1}(a_{n}b_{k}+b_{n}a_{k})c_{k}.
\end{equation*}
Now we can write $2\Im (ba)=2\sum_{i}e_{i}(b_{j}a_{k}-b_{k}a_{j})$ where
$(i,j,k)$ is a cyclic permutation of a quaternionic triple. Commuting
this with $c$ (which is equivalent to the action of $2\Im (ba)I$ on an element in a
matrix $H_{2}(\mk _{1})$) gives
\begin{equation*}
	2\Im (ba)=2(b_{j}a_{k}-a_{j}b_{k})c_{m}(e_{i}e_{m}-e_{m}e_{i})
\end{equation*}
Now $e_{i}e_{m}=-e_{m}e_{i}$, thus
\begin{equation*}
	[2\Im (ba),c]=\sum_{j}4((b_{j}a_{k}-b_{k}a_{j})c_{m}e_{j}
\end{equation*}
where $(e_{i},e_{m},e_{j})$ form a quaternionic triple. Therefore
$4(\mathbf{a}_{1}\mathbf{a}_{2}^{t}+\mathbf{a}_{2}\mathbf{a}_{1}^{t})$
is equivalent to $2\Im
(a_{2}a_{1})I$.
\item Now $\mathbf{b}^{t}\mathbf{a}=\sum_{i=1}^{\nu
-1}b_{i}a_{i}$. But
\begin{align*}
	(a,b)&=\Re (a\bar{b}) \\
	&= -\sum_{i=1}^{\nu -1}a_{i}b_{i}(e_{i})^{2} \\
	&= \sum_{i=1}^{\nu -1}a_{i}b_{i}
\end{align*}
as required. Clearly $\Re (ab)=-(a,b)$.
\end{enumerate} 
\end{proof}
\end{lemma}
Using these results and the fact that $G^{2}=I_{\nu -1}$ the two tables
are clearly equivalent and $\psi$ is a Lie algebra isomorphism.
\end{proof}

\subsection{Proof of equation~(\ref{eqn:parsnip3})}
We start by defining the structure and conformal algebras as they apply
to the $2\times 2$ hermitian matrix case. We will see that the structure
algebra is a subalgebra of the conformal algebra and the derivation
algebra a subalgebra of the structure algebra. This means that our
isomorphism for the conformal algebra includes the isomorphisms for the
structure algebra, and trivially for the derivation algebra.
From equations~(\ref{eqn:Str}) and~(\ref{eqn:der}) we can deduce that 
\begin{equation}
\label{eqn:str}
	\Str ^{\prime} H_{2}(\mk ) = H_{2}^{\prime }(\mk ) \ds A_{2}^{\prime }(\mk
) \ds \so (\mk ^{\prime })
\end{equation}
which is a Lie algebra with brackets defined by the statement that
$\Der (H_{2}(\mk ))  = A_{2}^{\prime }(\mk
) \ds \so(\mk ^{\prime })$ is a subalgebra. Denoting the elements of
the algebra by
\begin{align*}
	D &\in \Der (H_{2}(\mk )) \\
	H &\in H_{2}^{\prime }(\mk ) 
\end{align*}
we then have the brackets
\begin{align*}
	[D,H] &= D(H)\\
	[H,H^{\prime}] &= [L_{H},L_{H^{\prime}}]
\end{align*} 
where we consider the derivation $[L_{H},L_{H^{\prime}}]$ to be an element of
$\Der \mk $.
For the conformal algebra given by~\eqref{eqn:Con} we define the Lie
brackets by taking $\Str  \mj$  as a subalgebra and the brackets as 
below. We take the elements of $\Con  H_{2}(\mk )$ to be
\begin{align*}
	T &\in \Str \mj \\
	(X,Y) &\in [H_{2}(\mk )]^{2}.
\end{align*}
Then if $R \rightarrow R^{\ast }$ is an involutive automorphism which is
the identity on $\Der \mathbb{J}$ and multiplies elements of
$\mathbb{J}^{2}$ by $-1$ then the Lie brackets are
\begin{align*}
	[T,(X,Y)] &= (TX, T^{\ast }Y) \\
	[(X,0),(Y,0)] = [(0,X),(0,Y)] &= 0 \\
	[(X,0),(0,Y)] &= 2L_{XY} + 2[L_{X},L_{Y}]. 
\end{align*}
In the case of the conformal algebra for $2\times 2$ hermitian matrices
we substitute $H_{2}(\mk)$ for $\mj$ in equation~\eqref{eqn:Con} to obtain
\begin{equation*}
	\Con H_{2}(\mk) = \Str H_{2}(\mk) \ds [H_{2}(\mk )]^{2}
\end{equation*}
which we can expand to
\begin{equation}
\label{eqn:con }
	\Con  H_{2}(\mk ) =  H_{2}^{\prime }(\mk ) \ds A_{2}^{\prime }(\mk
) \ds \so (\mk ^{\prime }) \ds \mr \ds [H_{2}(\mk )]^{2}.
\end{equation}
We can now define the Lie brackets explicitly for $\Con H_{2}(\mk)$ by
considering the following elements
\begin{align*}
	  D &\in \Der H_{2}(\mk)  \\
	H &\in H_{2}^{\prime}(\mk) \\
	r,r^{\prime}&\in \mr  \\
	(X,Y)&\in [H_{2}(\mk )]^{2}. 
\end{align*}
Then the brackets are
\begin{align*}
	[D,r]=[r,H]&=[r,r^{\prime}]=0 \\
	[D,H]&=D(H) \\
	[D,(X,Y)]&=(D(X),D(Y)) \\
	[r,(X,Y)]&=(rX,-rY) \\
	[H,(X,Y)]&=(H\cdot X,-H\cdot Y)
\end{align*}
with the brackets for $(X,Y)$ defined as above.
We can also think of $\Str H_{2}(\mk)$ and $\Con H_{2}(\mk)$ in terms of
$2\times 2$ matrices over $H_{2}(\mk)$. Writing $\Str H_{2}(\mk)$ and $\Con
H_{2}(\mk)$ in this way gives 
\begin{align*}
	\Str H_{2}(\mk ) &= \Der H_{2}(\mk) \ds \left\{ \begin{pmatrix}
A & 0 \\ 0 & -A \end{pmatrix} \quad A \in H_{2}(\mk) \right\} \\
	\Con  H_{2}(\mk ) &=  \Der H_{2}(\mk) \ds \left\{ \begin{pmatrix}
A & B \\ C & -A \end{pmatrix} \quad A,B,C \in H_{2}(\mk) \right\}. 
\end{align*}

In this section we denote the elements of $L_2(\mk ,\mhs)$ by
\begin{align*}
	A&\in A_{2}^{\prime}(\mk ) \\
	S&\in \so (\mk ^{\prime}) \\
	B\otimes \tilde{i},C\otimes \tilde{j}, D\otimes \tilde{k} \in
H_{2}^{\prime}(\mk )\otimes \mhs ^{\prime} 
\end{align*}
along with the basis elements of $\so (\mhs ^{\prime})\cong \so (2,1)$ which we will call
$s_{12}, s_{13}$ and $s_{23}$, where
\begin{align*}
	s_{12}&=\begin{pmatrix} 0 & -1 & 0 \\ 1 & 0 & 0 \\ 0 & 0 & 0
\end{pmatrix} \\
	s_{13}&=\begin{pmatrix} 0 & 0 & -1 \\ 0 & 0 & 0 \\ -1 & 0 & 0
\end{pmatrix} \\
	s_{23}&=\begin{pmatrix} 0 & 0 & 0 \\ 0 & 0 & -1 \\ 0 & -1 & 0
\end{pmatrix}. 
\end{align*}
Now we define an isomorphism $\psi :L_2(\mk ,\mhs)\rightarrow \Con  H_{2}(\mk
)$ by
\begin{alignat*}
	 \psi \psi (A) &= A  & \quad	\psi (S) &= S 	 \\
	\psi (B\otimes \tilde{i}) &=B     \\
	\psi (C\otimes \tilde{j}) &= \tfrac{1}{2} (C,C), &\quad  	\psi (D\otimes \tilde{k}) &= \tfrac{1}{2}
(D,-D) \\
	\psi (s_{12}) &= \tfrac{1}{2} (I,-I), &\quad	\psi (s_{13}) &= \tfrac{1}{2}
(I,I) \\
	\psi (s_{23}) &= 1. 
\end{alignat*}
We can also define this in terms of our $4\times 4$ matrices by 
\begin{alignat*}
	 \psi \psi (A) &= A  & \quad	\psi (S) &= S 	 \\
	\psi (B\otimes \tilde{i}) &=\begin{pmatrix} B & 0 \\ 0 & -B \end{pmatrix}     \\
	\psi (C\otimes \tilde{j}) &= \begin{pmatrix} 0 & \tfrac{1}{2}C \\ 0 & 0 \end{pmatrix}, &\quad 
	\psi (D\otimes \tilde{k}) &= \begin{pmatrix} 0 & 0 \\ \tfrac{1}{2}D & 0 \end{pmatrix}
 \\
	\psi (s_{12}) &= \begin{pmatrix} 0 & 0 \\ \tfrac{1}{2}I & 0 \end{pmatrix},
&\quad 	\psi (s_{13}) &= \begin{pmatrix} 0 & \tfrac{1}{2}I \\ 0 & 0 \end{pmatrix}
 \\
	\psi (s_{23}) &= \begin{pmatrix} I & 0 \\ 0 & -I \end{pmatrix} . 
\end{alignat*}

The proof is a series of routine calculations of each product in both
algebras showing that the isomorphism holds in all cases,which can be
found on the next page.

Now $L_{2}(\mk, \mcs)$ is embedded in $L_{2}(\mk, \mhs)$ because $\mcs$
is embedded in $\mhs$ and $\psi$ maps $L_{2}(\mk, \mcs)$ to $\Str
^{\prime} H_{2}(\mk)$. Thus if we define \mbox{$\psi_{1}=\psi \mid L_{2}(\mk,
\mcs)$} then $\psi_{1}$ is an isomorphism between $L_{2}(\mk, \mcs)$ and $\Str
^{\prime} H_{2}(\mk)$. Similarly we can define $\psi_{2}=\psi \mid L_{2}(\mk,
\mr)$ as an isomorphism between $L_{2}(\mk, \mr)$ and $\Der H_{2}(\mk)$
by restricting $\psi$ to $L_{2}(\mk, \mr)$. Thus we have obtained proofs
for all the $2\times 2$ isomorphisms.

\begin{landscape}

\small
\begin{tabular}{|p{2cm}||p{2cm}|p{1.75cm}|p{2.25cm}|p{2cm}|p{2.25cm}|p{1.5cm}|p{1.5cm}|p{1.5cm}|}\hline 
  & $A^{\prime}$ & $D^{\prime}$ & $B^{\prime}\otimes \tilde{i}$ & $C^{\prime}\otimes \tilde{j}$
& $E^{\prime}\otimes \tilde{k}$ & $s_{12}$ & $s_{13}$ &
$s_{23}$ \\
 \hline \hline
  $A $ & $AA^{\prime}-A^{\prime}A$ & $-D^{\prime}(A)$ &
$(AB^{\prime}-B^{\prime}A)$ & $(AC^{\prime}-C^{\prime}A)$
& $(AE^{\prime}-E^{\prime}A)$ & $0$ & $0$ & $0$ \\
 & & & $\hspace{1cm} \otimes \tilde{i}$ & $\hspace{1cm} \otimes \tilde{j}$ & $\hspace{1cm} \otimes \tilde{k}$ &
& & \\
\hline
$D$ & $D(A^{\prime})$ & $DD^{\prime}-D^{\prime}D$ &
$D(B^{\prime})\otimes \tilde{i}$ & $D(C^{\prime})\otimes \tilde{j}$ & $D(E^{\prime})\otimes
\tilde{k}$ & $0$ & $0$ & $0$ \\
	\hline
$B\otimes \tilde{i}$ & $(BA^{\prime}-A^{\prime}B)$ &
$-D^{\prime }(B)\otimes \tilde{i}$ & $[L_{B},L_{B^{\prime }}]$ &
$ \frac{1}{2}\langle B,C^{\prime }\rangle s_{12} + $ & $ \frac{1}{2}\langle B,E^{\prime }\rangle s_{13} + $ & $B\otimes \tilde{j}$ & $B\otimes \tilde{k}$ &
$0$ \\
 & $\otimes \tilde{i}$ &  & & $(B\ast C^{\prime
})\otimes \tilde{k}$ & $(B\ast E^{\prime
})\otimes \tilde{j}$ & & & \\
	\hline
$C\otimes \tilde{j}$ & $(CA^{\prime}-A^{\prime}C)$ &
$-D^{\prime }(C)\otimes \tilde{j}$ &
$-\frac{1}{2}\langle C,B^{\prime }\rangle s_{12}-$ & $[L_{C},L_{C^{\prime }}]$ &
 $- \frac{1}{2}\langle C,E^{\prime }\rangle s_{23}-$ & $-C\otimes \tilde{i}$ & $0$ & $C\otimes \tilde{k}$ 
\\
 & $\otimes \tilde{j}$ & & $(C\ast B^{\prime
})\otimes \tilde{k}$ & & $(C\ast E^{\prime
})\otimes \tilde{i}$ & & & \\
	\hline
$E\otimes \tilde{k}$ & $(EA^{\prime}-A^{\prime}E)$ &
$-D^{\prime }(E)\otimes \tilde{k}$ &
$- \frac{1}{2}\langle E,B^{\prime }\rangle s_{13}-$ &
 $ \frac{1}{2}\langle E,C^{\prime }\rangle s_{23}+$ & $[L_{E},L_{E^{\prime }}]$ & $0$ &  
$E\otimes \tilde{i}$ & $E\otimes \tilde{j}$  \\
 & $\otimes \tilde{k}$ & & $(E\ast B^{\prime
})\otimes \tilde{j}$ & $(E\ast C^{\prime
})\otimes \tilde{i}$ & & & & \\
	\hline
$s_{12}$ & $0$ & $0$ & $-B^{\prime }\otimes \tilde{j}$ &
$C^{\prime }\otimes \tilde{i}$ & $0$ & $0$ & $s_{23}$ & $-s_{13}$
\\
	\hline
$s_{13}$ & $0$ & $0$ & $-B^{\prime }\otimes \tilde{k}$ & $0$ &
$-E^{\prime }\otimes \tilde{i}$ & $-s_{23}$ & $0$ & $-s_{12}$
\\
	\hline
$s_{23}$ & $0$ & $0$ & $0$ & $-C^{\prime }\otimes \tilde{k}$ &
$-E^{\prime }\otimes \tilde{j}$ & $s_{13}$ & $s_{12}$ & $0$
\\
	\hline  
    \end{tabular} \\
\vspace{0.5cm}
\begin{center}
The multiplication table for $L_{2}(\mk ,\mhs )$. \\
\end{center} 
\end{landscape}

\newpage
\begin{landscape}
\tiny
\begin{tabular}{|p{1.25cm}||p{2.145cm}|p{2.28cm}|p{2.5cm}|p{2.6cm}|p{2.5cm}|p{1.25cm}|p{1.25cm}|p{1.5cm}|}
\hline
   & $A^{\prime}$ & $D^{\prime}$ & $B^{\prime}$ & $\frac{1}{2}(C^{\prime},C^{\prime})$
& $\frac{1}{2}(E^{\prime},-E^{\prime})$ & $\frac{1}{2}(I,-I)$ & $\frac{1}{2}(I,I)$ &
$1$ \\
 \hline \hline
  $A $ & $AA^{\prime}-A^{\prime}A$ & $-D^{\prime}(A)$ &
$AB^{\prime}-B^{\prime}A$ & $\frac{1}{2}(AC^{\prime}-C^{\prime}A, $
& $\frac{1}{2}(AE^{\prime}-E^{\prime}A, $ & $0$ & $0$ & $0$
\\
 & & & & $AC^{\prime}-C^{\prime}A)$ & $E^{\prime}A-AE^{\prime})$ & & & \\
  \hline
$D$ & $D(A^{\prime})$ & $DD^{\prime}- D^{\prime}D$ &
$D(B^{\prime})$ & $\frac{1}{2}(D(C^{\prime}), D(C^{\prime}))$ &
$\frac{1}{2}(D(E^{\prime}), -D(E^{\prime}))$ & $0$ & $0$ & $0$ \\
	\hline
$B$ & $BA^{\prime}-A^{\prime}B$ &
$-D^{\prime }(B)$ & $[L_{B},L_{B^{\prime }}]$ &
$\frac{1}{4} \langle B,C^{\prime }\rangle (I, -I)+$ & $ \frac{1}{4}\langle B,E^{\prime }\rangle (I, I)+$ & $\frac{1}{2}(B,B)$ & $\frac{1}{2}(B,-B)$ & $0$
\\
 & & & & $(B\ast C^{\prime
},B\ast C^{\prime })$ & $(B\ast E^{\prime
},B\ast E^{\prime })$ & & & \\
	\hline
$\frac{1}{2}(C,C)$ & $\frac{1}{2}(CA^{\prime}-A^{\prime}C,$ &
$-\frac{1}{2}(D^{\prime }(C),D^{\prime }(C))$ &
$-\frac{1}{4} \langle C,B^{\prime }\rangle (I, -I)-$ & $[L_{C},L_{C^{\prime }}]$ &
 $-\frac{1}{2}\langle C,E^{\prime }\rangle 1-C\ast E^{\prime }$ & $-C$ & $0$ & $\frac{1}{2}(C,-C)$  \\
 & $CA^{\prime}-A^{\prime}C)$ & & $(C\ast B^{\prime
},C\ast B^{\prime })$ & & & & & \\
	\hline
$\frac{1}{2}(E,-E)$ & $\frac{1}{2}(EA^{\prime}-A^{\prime}E,$ &
$-\frac{1}{2}(D^{\prime }(E),D^{\prime }(E))$ &
$-\frac{1}{4} \langle E,B^{\prime }\rangle (I, I)-$ &
 $\frac{1}{2}\langle E,C^{\prime }\rangle 1+E\ast C^{\prime }$ & $[L_{E},L_{E^{\prime }}]$ & $0$ &  
$E$ & $\frac{1}{2}(E,E)$  \\
 & $A^{\prime}E-EA^{\prime})$ & & $(E\ast B^{\prime },E\ast B^{\prime
})$ & & & & & \\
	\hline
$\frac{1}{2}(I,-I)$ & $0$ & $0$ & $-\frac{1}{2}(B^{\prime },B^{\prime })$ &
$C^{\prime }$ & $0$ & $0$ & $1$ & $-\frac{1}{2}(I,I)$ \\
	\hline
$\frac{1}{2}(I,I)$ & $0$ & $0$ & $-\frac{1}{2}(B^{\prime },-B^{\prime })$ & $0$ &
$-E^{\prime }$ & $-1$ & $0$ & $-\frac{1}{2}(I,-I)$ \\
	\hline
$1$ & $0$ & $0$ & $0$ & $-\frac{1}{2}(C^{\prime },-C^{\prime })$ &
$-\frac{1}{2}(E^{\prime },E^{\prime })$ & $\frac{1}{2}(I,I)$ & $\frac{1}{2}(I,-I)$ & $0$ \\
	\hline
      \end{tabular} \\
\vspace{0.5cm}
\begin{center}
The multiplication table for $\Con H_{2}(\mk )$. 
\end{center}

\end{landscape}

\section{Magic squares of $3\times 3$ Matrix Algebras: Proofs}
\label{sec:proofs3}
In this section we extend the results of the last section to the
$3\times 3$ matrix case. We then develop these ideas by showing the maximal
compact subalgebras for each of the exceptional Lie algebras that appear
in the magic square.

We begin by showing that 
\begin{equation}
\label{eqn:carrot}
L_{3}(\mk ,\mhs )\cong \Con H_{3}(\mk )
\end{equation}
and then show the maximal compact subalgebras are as stated previously.

\subsection{Proof of equation~(\ref{eqn:carrot})}
In the case of $3\times 3$ matrices we know from~(\ref{eqn:Str}) and
Theorem~\ref{theorem:apple} that
\begin{equation*}
	\Str ^{\prime }H_{3}(\mk )=L_{3}^{\prime}(\mk )\ds\Der \mk
\end{equation*}
and also from~(\ref{eqn:Con}), 
\begin{equation*}
	\Con H_{3}(\mk )=\Str H_{3}(\mk )\ds [H_{3}(\mk )]^{2}.
\end{equation*}
Thus 
\begin{equation*}
	\Con H_{3}(\mk ) = \Der \mk \ds A_{3}^{\prime }(\mk ) \ds
H_{3}^{\prime }(\mk ) \ds \mr \ds [H_{3}(\mk )]^{2}.
\end{equation*}
We have
\begin{equation*}
	L_{3}(\mk ,\mhs )=\Der H_{3}(\mk ) \ds H_{3}^{\prime }(\mk
)\otimes \mhs ^{\prime} \ds \Der \mhs.
\end{equation*}
It is known from~\cite{Sudbery84} that $\Der \mh = C(\mh ^{\prime})$ and
can be shown similarly that $\Der \mhs = C(\mhs ^{\prime})$.
Thus, using Theorem~\ref{theorem:apple}
\begin{equation*}
	L_{3}(\mk ,\mhs )=\Der \mk \ds A_{3}^{\prime}(\mk ) \ds H_{3}^{\prime }(\mk
)\otimes \mhs ^{\prime} \ds C(\mhs ^{\prime }).
\end{equation*}
If we consider the elements of $L_{3}(\mk ,\mhs )$ 
\begin{align*}
	A&\in A_{3}^{\prime }(\mk ) \\
	D&\in \Der (\mk ) \\
	B\otimes \tilde{i},C\otimes \tilde{j}, E\otimes \tilde{k} &\in
H_{3}^{\prime}(\mk )\otimes \mhs ^{\prime} \\
	C_{\tilde{i}}, C_{\tilde{j}},C_{\tilde{k}} &\in C(\mhs ^{\prime}) 
\end{align*}
and the elements of $\Con H_{3}(\mk )$ 
\begin{align*}
	A &\in A_{3}^{\prime }(\mk ) \\
	B &\in H_{3}^{\prime }(\mk ) \\
	D &\in \Der (\mk ) \\
	r &\in \mr \\
	(X,Y) &\in [H_{3}(\mk )]^{2}
\end{align*}
then we can define an isomorphism $\phi : L_{3}(\mk ,\mhs)\to \Con
H_{3}(\mk )$ by :
\begin{alignat*}
	\phi \phi (A) &= A  & \quad	\phi (D) &= D 	 \\
	\phi (B\otimes \tilde{i}) &=B     \\
	\phi (C\otimes \tilde{j}) &= \tfrac{1}{2} (C,C) &\quad  	\phi (E\otimes \tilde{k})
&= \tfrac{1}{2}(E,-E) \\
	\phi (C_{\tilde{k}}) &= (I,-I) &\quad	\phi (C_{\tilde{j}}) &= (I,I) \\
	\phi (C_{\tilde{i}}) &= 2. 
\end{alignat*}
Again we can express $\Con H_{3}(\mk)$ and $\Str ^{\prime} H_{3}(\mk )$ in terms
of $2\times 2$ matrices over $H_{3}(\mk )$ given by
\begin{align*}
	\Str ^{\prime}H_{3}(\mk ) &= \Der H_{3}(\mk) \ds \left\{ \begin{pmatrix}
A & 0 \\ 0 & -A \end{pmatrix} \quad A \in H_{3}^{\prime}(\mk) \right\} \\
	\Con  H_{3}(\mk ) &=  \Der H_{3}(\mk) \ds \left\{ \begin{pmatrix}
A & B \\ C & -A \end{pmatrix} \quad A,B,C \in H_{3}(\mk) \right\}. 
\end{align*}
Further we can also define our isomorphism in terms of these matrices
\begin{alignat*}
	\phi \phi (A) &= A  & \quad	\phi (D) &= D 	 \\
	\phi (B\otimes \tilde{i}) &=\begin{pmatrix} B & 0 \\ 0 & -B \end{pmatrix}     \\
	\phi (C\otimes \tilde{j}) &=\begin{pmatrix} 0 & \tfrac{1}{2}C \\
0 & 0 \end{pmatrix} &\quad  	\phi (E\otimes \tilde{k})
&= \begin{pmatrix} 0 & 0 \\ \tfrac{1}{2}E & 0 \end{pmatrix} \\
	\phi (C_{\tilde{k}}) &= \begin{pmatrix} 0 & 0 \\ I & 0 \end{pmatrix}
&\quad	\phi (C_{\tilde{j}}) &= \begin{pmatrix} 0 & I \\ 0 & 0 \end{pmatrix} \\
	\phi (C_{\tilde{i}}) &= \begin{pmatrix} 2I & 0 \\ 0 & -2I \end{pmatrix}. 
\end{alignat*}

We can again show this is a Lie algebra isomorphism by comparisons of
the multiplication tables on the following page for $L_{3}(\mk ,\mhs)$ and $\Con
H_{3}(\mk )$.
\newpage
\begin{landscape}
\tiny

\begin{center}

\begin{tabular}{|p{1cm}||p{2cm}|p{1.75cm}|p{2cm}|p{2cm}|p{2cm}|p{1.5cm}|p{1.5cm}|p{1.5cm}|}
\hline
   & $A^{\prime}$ & $D^{\prime}$ & $B^{\prime}\otimes \tilde{i}$ & $C^{\prime}\otimes \tilde{j}$
& $E^{\prime}\otimes \tilde{k}$ & $C_{\tilde{k}}$ & $C_{\tilde{j}}$ &
$C_{\tilde{i}}$ \\
 \hline \hline
  $A $ & $AA^{\prime}-A^{\prime}A$ & $-D^{\prime}(A)$ &
$(AB^{\prime}-B^{\prime}A)\otimes \tilde{i}$ & $(AC^{\prime}-C^{\prime}A)\otimes \tilde{j}$
& $(AE^{\prime}-E^{\prime}A)\otimes \tilde{k}$ & $0$ & $0$ & $0$ \\
  \hline
$D$ & $D(A^{\prime})$ & $DD^{\prime}-D^{\prime}D$ &
$D(B^{\prime})\otimes \tilde{i}$ & $D(C^{\prime})\otimes \tilde{j}$ & $D(E^{\prime})\otimes
\tilde{k}$ & $0$ & $0$ & $0$ \\
	\hline
$B\otimes \tilde{i}$ & $(BA^{\prime}-A^{\prime}B)\otimes \tilde{i}$ &
$-D^{\prime }(B)\otimes \tilde{i}$ & $[L_{B},L_{B^{\prime }}]$ &
$\frac{1}{6} \langle B,C^{\prime }\rangle C_{\tilde{k}}+(B\ast C^{\prime
})\otimes \tilde{k}$ & $\frac{1}{6} \langle B,E^{\prime }\rangle C_{\tilde{j}}+(B\ast E^{\prime
})\otimes \tilde{j}$ & $2B\otimes \tilde{j}$ & $2B\otimes \tilde{k}$ &
$0$ \\
	\hline
$C\otimes \tilde{j}$ & $(CA^{\prime}-A^{\prime}C)\otimes \tilde{j}$ &
$-D^{\prime }(C)\otimes \tilde{j}$ &
$-\frac{1}{6} \langle C,B^{\prime }\rangle C_{\tilde{k}}-(C\ast B^{\prime
})\otimes \tilde{k}$ & $[L_{C},L_{C^{\prime }}]$ &
 $-\frac{1}{6} \langle C,E^{\prime }\rangle C_{\tilde{i}}-(C\ast E^{\prime
})\otimes \tilde{i}$ & $-2C\otimes \tilde{i}$ & $0$ & $2C\otimes \tilde{k}$  \\
	\hline
$E\otimes \tilde{k}$ & $(EA^{\prime}-A^{\prime}E)\otimes \tilde{k}$ &
$-D^{\prime }(E)\otimes \tilde{k}$ &
$-\frac{1}{6} \langle E,B^{\prime }\rangle C_{\tilde{j}}-(E\ast B^{\prime
})\otimes \tilde{j}$ &
 $\frac{1}{6} \langle E,C^{\prime }\rangle C_{\tilde{i}}+(E\ast C^{\prime
})\otimes \tilde{i}$ & $[L_{E},L_{E^{\prime }}]$ & $0$ &  
$2E\otimes \tilde{i}$ & $2E\otimes \tilde{j}$  \\
	\hline
$C_{\tilde{k}}$ & $0$ & $0$ & $-2B^{\prime }\otimes \tilde{j}$ &
$2C^{\prime }\otimes \tilde{i}$ & $0$ & $0$ & $2C_{\tilde{i}}$ & $-2C_{\tilde{j}}$
\\
	\hline
$C_{\tilde{j}}$ & $0$ & $0$ & $-2B^{\prime }\otimes \tilde{k}$ & $0$ &
$-2E^{\prime }\otimes \tilde{i}$ & $-2C_{\tilde{i}}$ & $0$ & $-2C_{\tilde{k}}$
\\
	\hline
$C_{\tilde{i}}$ & $0$ & $0$ & $0$ & $-2C^{\prime }\otimes \tilde{k}$ &
$-2E^{\prime }\otimes \tilde{j}$ & $2C_{\tilde{j}}$ & $2C_{\tilde{k}}$ & $0$
\\
	\hline
      \end{tabular}\\
\vspace{0.5cm}
The multiplication table for $L_{3}(\mk ,\mhs )$. \\
\end{center}

\end{landscape}
\begin{landscape}
\tiny
\begin{center}

\begin{tabular}{|p{1.25cm}||p{2.4cm}|p{2.28cm}|p{2.5cm}|p{2.6cm}|p{2.6cm}|p{1.25cm}|p{1.25cm}|p{1.5cm}|}
\hline
   & $A^{\prime}$ & $D^{\prime}$ & $B^{\prime}$ & $\frac{1}{2}(C^{\prime},C^{\prime})$
& $\frac{1}{2}(E^{\prime},-E^{\prime})$ & $(I,-I)$ & $(I,I)$ &
$2$ \\
 \hline \hline
  $A $ & $AA^{\prime}-A^{\prime}A$ & $-D^{\prime}(A)$ &
$AB^{\prime}-B^{\prime}A$ & $\frac{1}{2}(AC^{\prime}-C^{\prime}A,$
& $(AE^{\prime}-E^{\prime}A,$ & $0$ & $0$ & $0$
\\
 & & & & $AC^{\prime}-C^{\prime}A)$ & $-AE^{\prime}+E^{\prime}A)$ & & & \\
  \hline
$D$ & $D(A^{\prime})$ & $DD^{\prime}- D^{\prime}D$ &
$D(B^{\prime})$ & $\frac{1}{2}(D(C^{\prime}), D(C^{\prime}))$ &
$\frac{1}{2}(D(E^{\prime}), -D(E^{\prime}))$ & $0$ & $0$ & $0$ \\
	\hline
$B$ & $BA^{\prime}-A^{\prime}B$ &
$-D^{\prime }(B)$ & $[L_{B},L_{B^{\prime }}]B$ &
$\frac{1}{6} \langle B,C^{\prime }\rangle (I, -I)+$ & $\frac{1}{6} \langle B,E^{\prime }\rangle (I, I)+$ & $(B,B)$ & $(B,-B)$ & $0$ \\
 & & & & $(B\ast C^{\prime
},B\ast C^{\prime })$ & $(B\ast E^{\prime
},B\ast E^{\prime })$ & & & \\
	\hline
$\frac{1}{2}(C,C)$ & $\frac{1}{2}(CA^{\prime}-A^{\prime}C,$ &
$-\frac{1}{2}(D^{\prime }(C),D^{\prime }(C))$ &
$-\frac{1}{6} \langle C,B^{\prime }\rangle (I, -I)-$ & $[L_{C},L_{C^{\prime }}]$ &
 $-\frac{1}{3} \langle C,E^{\prime }\rangle -C\ast E^{\prime }$ & $-2C$ & $0$ & $(C,-C)$ 
\\
 & $CA^{\prime}-A^{\prime}C)$ & & $(C\ast B^{\prime
},C\ast B^{\prime })$ & & & & & \\
	\hline
$\frac{1}{2}(E,-E)$ & $\frac{1}{2}(EA^{\prime}-A^{\prime}E,$ &
$-\frac{1}{2}(D^{\prime }(E),D^{\prime }(E))$ &
$-\frac{1}{6} \langle E,B^{\prime }\rangle (I, I)-$ &
 $\frac{1}{3} \langle E,C^{\prime }\rangle +E\ast C^{\prime }$ & $[L_{E},L_{E^{\prime }}]$
&  $0$ &  $2E$ & $(E,E)$  \\
 & $-EA^{\prime}+A^{\prime}E)$ & & $(E\ast B^{\prime },E\ast B^{\prime
})$ & & & & & \\
	\hline
$(I,-I)$ & $0$ & $0$ & $-(B^{\prime },B^{\prime })$ &
$2C^{\prime }$ & $0$ & $0$ & $4$ & $-2(I,I)$ \\
	\hline
$(I,I)$ & $0$ & $0$ & $-(B^{\prime },B^{\prime })$ & $0$ &
$-2E^{\prime }$ & $-4$ & $0$ & $-2(I,-I)$ \\
	\hline
$2$ & $0$ & $0$ & $0$ & $-(C^{\prime },-C^{\prime })$ &
$-(E^{\prime },E^{\prime })$ & $2(I,I)$ & $2(I,-I)$ & $0$ \\
	\hline
    \end{tabular} \\
\vspace{0.5cm}
  The multiplication table for $\Con H_{3}(\mk )$.
    \end{center}
\end{landscape}

Invoking the same method as used in the similar proof for $2\times 2$
matrices we note that we can define  $\phi _1 : L_{3}(\mk ,\mcs)\to \Str
^{\prime } H_{3}(\mk )$ by:
\begin{alignat*}
	\phi \phi _{1}(A) &= A,  & \quad	\phi _{1}(D) &= D 	 \\
	\phi _{1}(B\otimes \tilde{i}) &=B     
\end{alignat*}
and also, trivially, $\phi _{2}: L_{3}(\mk ,\mr)\to \Der H_{3}(\mk )$ by :
\begin{equation*}
	\phi _{2}(A) = A,   \quad	\phi _{2} (D) = D 	 
\end{equation*}
Thus we have proved the following.
\begin{theorem} 
\begin{align*}
	L_3(\mk ,\mr ) &\cong \Der H_3 (\mk ) \\
	L_3(\mk ,\mcs) &\cong \Str ^{\prime }H_3(\mk )\\
	L_3(\mk ,\mhs) &\cong \Con  H_3(\mk ) 
\end{align*}
hold for $\mk = \mr ,\mc ,\mh $ and $\mo $.
\end{theorem}

The relation to the matrix algebras described in the Introduction in
equation~(\ref{eqn:iso2}) is:
\begin{theorem}
\begin{align*}
 L_{3}(\mk ,\mr)&\cong \su(3,\mk) \\
 L_{3}(\mk ,\mcs)&\cong \sl(3,\mk)  \\
 L_{3}(\mk ,\mhs)&\cong \sp(6,\mk).
\end{align*}
\end{theorem}
Thus we have the table
\begin{center}
\vspace{0.5cm} \begin{tabular}{|c||c|c|c|c|}
\hline
     & $\mr$ & $\mc$  & $\mh$ & $\mo$ \\
 \hline \hline
  $\Der H_3(\mk )\cong L_3(\mk ,\mr)$ & $\su(3,\mr)$ & $\su(3,\mc)$ &
$\su(3,\mh )$ & $\su(3,\mo)$ \\
  \hline
   $\Str ^\prime  H_3(\mk )\cong L_3(\mk ,\mcs ) $ & $\sl(3,\mr )$ &
$\sl(3,\mc )$ & $\sl(3,\mh )$ & $\sl(3,\mo)$ \\
   \hline
    $\Con  H_3(\mk )\cong L_3(\mk ,\mhs ) $ & $\sp(6,\mr )$ & $\sp(6,\mc)$ &
$\sp(6,\mh )$ & $\sp(6,\mo)$ \\
    \hline
    \end{tabular} \vspace{0.5cm}
    \end{center}
\begin{proof}
\begin{enumerate}
\item  We recall that $\su (3,\mk)$ is
the set of matrices with entries in $\mk$ satisfying $A^{\dag}G+GA=0$,
i.e. the set of anti-hermitian matrices. In this case we have
$L_{3}(\mk, \mr) = A_{3}^{\prime}(\mk )\ds \Der \mk $. Take $A\in A_{3}^{\prime}(\mk
)$ and $D \in \Der \mk$. Then $\psi : L_{3}(\mk ,\mr) \rightarrow \su
(3,\mk )$ by
\begin{equation*}
	\psi (A+D)=A+DI
\end{equation*}
where $I$ is the $3\times 3$ matrix identity. 
\item We define $\sl (3,\mk )$ to be the family of $3\times 3$ matrices
with entries in $\mk$ with the real part of the trace equal to zero.
Also
\begin{equation*}
	L_{3}(\mk,\mcs)=A_{3}^{\prime}(\mk) \ds \Der \mk \ds
H_{3}^{\prime}(\mk )\otimes \mcs .
\end{equation*}
If we take $A\in A_{3}^{\prime}(\mk)$, $D\in \Der \mk $ and $H \in
H_{3}^{\prime}(\mk)$ (since we can regard the tensor product $H_{3}^{\prime}(\mk )\otimes
\mcs$ as being one copy of $H_{3}^{\prime}(\mk)$) then the isomorphism
$\phi : L_{3}(\mk ,\mcs)\rightarrow \sl (3,\mk)$ can be written
\begin{equation*}
	\phi (A+D+H ) = A+DI+H.
\end{equation*} 
\item The Lie algebra $\sp (6,\mk )$ is defined to be the Lie algebra
of $6\times 6$ matrices satisfying the equation $A^{\dag}J+JA=0$ with
the added condition that the trace is also zero. This can be written as
the matrix $\begin{pmatrix} A & B \\ C & -A^{\dag} \end{pmatrix}$ where
$A$, $B$ and $C$ are $3\times 3$ block matrices and $B$ and $C$ are 
hermitian. We have
\begin{equation*}
	L_{3}(\mk ,\mhs)= A_{3}^{\prime}(\mh) \ds \Der \mh \ds
H_{3}^{\prime}(\mh) \otimes \mhs ^{\prime} \ds \Der \mhs .
\end{equation*}
Taking $A\in A_{3}^{\prime}(\mk)$, $D\in \Der \mk $, $H_{1} \otimes
\tilde{i}$, $H_{2} \otimes
\tilde{j}$ and $H_{3} \otimes
\tilde{k} \in H_{3}^{\prime}(\mk)\otimes \mhs$ and $r_{1}C_{\tilde{i}}+
r_{2}C_{\tilde{j}}+ r_{3}C_{\tilde{k}}\in \Der \mhs$ then the
isomorphism $\chi :L_{3}(\mk,\mhs)\rightarrow \sp(6,\mk)$ can be written
explicitly as
\begin{multline*}
	\chi (A+DI+H_{1}\otimes \tilde{i}+H_{2}\otimes \tilde{j}+H_{3}\otimes
\tilde{k}+ r_{1}C_{\tilde{i}}+
r_{2}C_{\tilde{j}}+ r_{3}C_{\tilde{k}})= \\
\begin{pmatrix} A+DI+H_{1}+\tfrac{1}{3}r_{1}I &
(H_{2}+r_{2}I-H_{3}-\tfrac{1}{3}r_{3}I) \\
(H_{2}+r_{2}I+H_{3}+\tfrac{1}{3}r_{3}I) & A+DI-H_{1}-\tfrac{1}{3}r_{1}I \end{pmatrix}.
\end{multline*}	
\end{enumerate}
\end{proof}

\subsection{Maximal Compact Subalgebras}
A semi-simple Lie algebra is called \emph{compact} if it has a
negative-definite killing form. It is called \emph{non-compact} if its killing
form is not negative-definite.

A non-compact real form, $\g$, of a semi-simple complex Lie algebra, $L$, has a
\emph{maximal compact subalgebra} $\f$ with an \emph{orthogonal complementary
subspace} $\p$ such that $\g = \f \ds \p $ and the brackets
\begin{align*}
	[\f ,\f ] & \subseteq \f \\
	[\f ,\p ] & \subseteq \p \\
	[\p ,\p ] & \subseteq \f \\
	(\f ,\p ) &= 0
\end{align*}
are satisfied. We denote by $(,)$ the killing
form of $L$. There exists an involutive automorphism $\sigma : \g
\rightarrow \g$ such that $\f $ and $\p $ are eigenspaces of $\sigma $ with
eigenvalues $+1$ and $-1$ respectively. A compact real form, $\g
^{\prime}$, of $L$ will also contain $\f$ as a compact subalgebra of $\g
^{\prime}$ but clearly in this case the maximal compact subalgebra will be $\g
^{\prime}$ itself. We can obtain $\g ^{\prime}$ from $\g$ by keeping the
same brackets in $[\f ,\f ]$ and $[\f ,\p ]$ but multiplying the
brackets in $[\p ,\p ]$ by $-1$, i.e. by performing the \emph{Weyl
unitary trick} (putting $\g
^{\prime}= \f \ds i\p $).

We will now give an overview of the method used to show that the
algebras given in the table on page~\pageref{maximal} are maximal
compact, which is essentially the same in each case. We know that
$L_{3}(\mk _{1},\mk _{2})$ gives a compact real form of each Lie algebra (from~\cite{JacobsonELA}).
Thus if  $L_{3}(\mk _{1},\mks _{2})$ shares a common subalgebra with $L_{3}(\mk _{1},\mk
_{2})$, say $\f $, where
\begin{align*}
	L_{3}(\mk _{1},\mk _{2}) &= \f \ds \p _{1} \\
	L_{3}(\mk _{1},\mks _{2})&= \f \ds \p _{2},
\end{align*}
and the brackets in $[\f , \p _{1}]$ are the same as those in $[\f ,\p
_{2}]$ but the brackets in $[\p _{1}, \p _{1}]$ are $-1$ times the
equivalent brackets in $[\p _{2},\p _{2}]$, then $\f $ will be the
maximal compact subalgebra of $L_{3}(\mk _{1},\mks _{2})$ and $\p _{2}$
will be its orthogonal complementary subspace. Moreover, because of the
nature of the split composition algebras, this sign change in the
brackets will reflect precisely the change in sign in the Cayley-Dickson
process when moving from the non-split to the split form of the
composition algebra. 

We will briefly consider the nature of $\Der \mo$ and $\Der \mos$ since
the structure of these algebras form a fundamental part of this proof.
It is well known that $\Der \mo \cong G_{2}$ (see~\cite{Sudbery84}). The
derivation algebra of $\mo$ is the Lie algebra of the automorphism group
of $\mo$. The derivations can be split into two types, those that are the
infinitesimal versions of the automorphisms of $\mo$ fixing the complex
subspace $\mc = \mr \ds i\mr $ and those which are the infinitesimal
versions of automorphisms fixing $\mr$. The derivations of the first
type (leaving $i$ invariant) form a subalgebra isomorphic to $\su (3)$.

Further we can express the elements of $\Der \mo$ in terms of pairs of basis elements $s_{ij}$ of $\so (7)$, as defined previously, since $\Der \mo \subset \so (7)$. Similarly
$\Der \mos \subset \so (4,3)$, giving a representation of $\Der \mos$ in
terms of pairs of basis elements of $\so (4,3)$.
It can be shown that the 14 elements of $\Der \mo $ (by a method similar to that found in ~\cite{Dixon}) are
\begin{align*}
	g_{1}&=s_{23}-s_{45} &\quad g_{2}&=s_{45}+s_{67} \\
	g_{3}&=s_{25}-s_{34} &\quad g_{4}&=-s_{27}-s_{36} \\
	g_{5}&=-s_{47}-s_{56} &\quad g_{6}&=-s_{24}-s_{35} \\
	g_{7}&=-s_{26}+s_{37} &\quad g_{8}&=-s_{46}+s_{57} \\
	g_{9}&=-s_{12}+s_{47} &\quad g_{10}&=s_{13}+s_{57} \\
	g_{11}&=s_{14}+s_{27} &\quad g_{12}&=-s_{15}+s_{37} \\
	g_{13}&=s_{16}+s_{25} &\quad g_{14}&=-s_{17}+s_{24}. 
\end{align*}
The derivations in $\Der \mos$ are obtained directly from these
by writing, for example
\begin{equation*}
	\tilde{g_{1}}=\tilde{s}_{23}-\tilde{s}_{45},
\end{equation*} 
where, if $S_{ij}$ is the matrix representation of $s_{ij}$ then
$G_{2}S_{ij}$ is the matrix representation of $\tilde{s_{ij}}$ (recall
that $G_{2}$ is the metric for $\mk _{2}$). The
other elements of $\Der \mos$ are obtained by a similar method. Then $\Der \mo $
and $\Der \mos$ have a common subalgebra $\so (3)\ds \so (3)$
(see, for example,~\cite{Gilmore}) which is the subalgebra with basis
elements $\{ g_{1}, g_{9}, g_{10}, g_{2}, g_{5}, g_{8} \}$, these being
invariant under multiplication by  the metric $G_{2}$. We will denote by
$\tilde{G}_{2}$ the form of the exceptional Lie algebra $G_{2}$
isomorphic to $\Der \mos$.

We now state explicity the result we are about to prove.
\begin{theorem}
The maximal compact subalgebras, stated explicity, are of the forms given
in the table below.
\begin{center}
\vspace{0.5cm} \begin{tabular}{|c|c|c|}
\hline
	$\g $ & $\f $ &  \\
\hline

	$E_{6,1}$ & $F_4$ & $\Der H_{3}(\mo )$ \\
\hline
	$E_{7,1}$ & $E_{6} \oplus \so (2)$ & $\Der H_{3}(\mo
)\ds H^{\prime }_{3}(\mo) \otimes \{ i\} \ds \{ C_{i} \} $  \\
\hline
	$E_{8,1}$ & $E_{7} \oplus \so (3)$ & $\Der H_{3}(\mo
)\ds H^{\prime }_{3}(\mo) \otimes \mh ^{\prime} \ds M_{G_{2}}$  \\
\hline 

	$E_{6,2}$ & $\sq (3) \oplus \so (3)$ & $\Der H_{3}(\mc )\ds
H_{3}^{\prime} (\mc) \otimes \mh ^{\prime} \ds M_{G_{2}}$ \\
\hline
	$E_{7,2}$ & $\su (6) \oplus \so (3)$ & $\Der H_{3}(\mh )\ds
H_{3}^{\prime} (\mh) \otimes \mh ^{\prime} \ds M_{G_{2}}$  \\
\hline
\end{tabular} \vspace{0.5cm}.
\end{center}
where $M_{G_{2}} (= \{ g_{1},g_{9},g_{10},g_{2},g_{5},g_{8} \} )$, is the
maximal compact subalgebra of $\tilde{G_{2}}$.

\begin{proof}
We consider first each of the $ _{,1}$ type algebras and then move on to
the $_{,2}$ type. Denote by $A^{N}$ the non- compact form of the algebra $A$ and by $A^{C}$ the compact form of $A$. 
\begin{enumerate}
\item  For $E_{6,1}^{N}$ 
\begin{align*}
	\f &= \Der H_{3}(\mo ) \\
	\p &= H_{3}^{\prime }(\mo ) \otimes \mcs ^{\prime}.
\end{align*}
$E_{6}^{C}$ also has $\f $ as a subalgebra but in this case the remaining subspace is $H_{3}^{\prime }(\mo ) \otimes \mcs ^{\prime}$.
 Thus there is only one set of brackets to
check. If we consider $H_{1}\otimes \tilde{i},H_{2}\otimes \tilde{i} \in H_{3}(\mo )\otimes
\mcs$ and $H_{1}\otimes i,H_{2}\otimes i \in H_{3}(\mo )\otimes
\mc$ then clearly, using the definitions for the brackets found on page
\pageref{brackets}
\begin{align*}
	[H_{1}\otimes \tilde{i},H_{2}\otimes \tilde{i}] &=
[L_{H_{1}},L_{H_{2}}] \\
	 [H_{1}\otimes i,H_{2}\otimes i] &=
-[L_{H_{1}},L_{H_{2}}] ,
\end{align*} as required. 
\item In the case of $E_{7,1}^{N}$ the orthogonal complementary subspace is
\begin{equation*}
	\tilde{\p } = H_{3}^{\prime}(\mo )\otimes \{
\tilde{j},\tilde{k} \} \ds \{ C_{\tilde{j}},C_{\tilde{k}} \} ,	
\end{equation*}
where $\f $ is the maximal compact subalgebra given in the table above. Then $\f $ is also a subalgebra in $E_{7}^{C}$ and the remaining subspace is $\p = H_{3}^{\prime}(\mo )\otimes \{ j,k\} \ds \{ C_{j},C_{k}
\}$. Now
\begin{align*}
	[C_{\tilde{j}},C_{\tilde{k}}] &= -2C_{i} &\quad [C_{j},C_{k}] &=
2C_{i} \\
	[C_{\tilde{j}},H_{1}\otimes \tilde{k}]&= -H_{1}\otimes 2i &\quad
[C_{j},H_{1}\otimes k]&= H_{1}\otimes 2i \\
	[C_{\tilde{k}},H_{1}\otimes \tilde{j}]&=H_{1}\otimes 2i &\quad
[C_{k},H_{1}\otimes j]&= -H_{1}\otimes 2i.
\end{align*}
Recall that
\begin{multline*}
\qquad [H_{1}\otimes x,H_{2} \otimes y] = 2\tr(H_{1}H_{2})D_{x,y} + \\(H\ast
G)\otimes \Im (xy) 
- Re(x\bar{y})[L_{H_{1}},L_{H_{2}}].
\end{multline*}
We have to consider the two cases (1) when $x=y$ and (2) when $x\perp y$.
Case (1) is considered in $E_{6}$. Case (2) gives
\begin{align*}
	[H_{1}\otimes \tilde{j},H_{2}\otimes \tilde{k}] &= -2\tr
(H_{1}H_{2})C_{i}-(H_{1}\ast H_{2})\otimes i \\
	[H_{1}\otimes j,H_{2}\otimes k] &= 2\tr
(H_{1}H_{2})C_{i}+(H_{1}\ast H_{2})\otimes i. 
\end{align*}
Clearly, $\p$ and $\tilde{\p }$ are orthogonal complementary subspaces
with their brackets with themselves giving opposite signs and thus we
can deduce that the choice of maximal compact subalgebra is correct. 
\item In $E_{8}^{N}$ we have the orthogonal complementary subspace
\begin{equation*}
	\tilde{\p } = H_{3}^{\prime}(\mo )\otimes \{ \tilde{l},\tilde{il},\tilde{jl},\tilde{kl} \} \ds  \{ \tilde{g_{a}}
\mid 
a=3,4,6,7,11,12,13,14 \} 
\end{equation*}
to the maximal compact subalgebra $\f $ as shown in the previous table. Then $\f $ is also a subalgebra of $E_{8}^{C}$ and the remaining subspace in $E_{8}^{C}$ will be $\p = H_{3}^{\prime}(\mo )\otimes \{ l,il,jl,kl \} \ds  \{ g_{a} \mid a=3,4,6,7,11,12,13,14 \}$. 
For convenience we will label the orthogonal complementary subspaces of
$G_{2}$ $\p_{G_{2}}$ and $\tilde{\p}_{G_{2}}$ for the compact and
non-compact cases respectively. The calculations for $E_{8}$ are much
the same as those for $E_{7}$. For brackets between $H_{3}^{\prime}(\mo )\otimes \{ l,il,jl,kl
\}$ and itself we again have two cases with $x=y$ and $x\perp y$. These
are resolved in the same way as before.
To calculate the brackets between $H_{3}^{\prime}(\mo )\otimes \{ l,il,jl,kl
\}$ and $\p_{G_{2}}$ and between $\p_{G_{2}}$ and itself involves a set
of long but relatively simple calculations involving $s_{ij}$ and
$\tilde{s_{ij}}$. These produce the signs as expected, however, in the
interest of the rainforests, we have not reproduced them here. 
\end{enumerate}

Now notice that these proofs do not in fact involve the matrices in
$H_{3}(\mk)$ since the change between compactness and non-compactness
does not involve $\mk_{1}$ but only $\mk_{2}$.  Thus since the
orthogonal compact subspaces of $E_{6,2}$ and $E_{7,2}$ are
\begin{align*}
	\tilde{\p_{1}} &= H_{3}^{\prime}(\mc )\otimes \{
\tilde{l},\tilde{il},\tilde{jl},\tilde{kl} \} \ds \tilde{\p}_{G_{2}}
\\
	\tilde{\p_{1}} &= H_{3}^{\prime}(\mh )\otimes \{
\tilde{l},\tilde{il},\tilde{jl},\tilde{kl} \} \ds \tilde{\p}_{G_{2}},
\end{align*}
the proofs for these maximal compact subalgebras are contained within
the that of $E_{8}$

Thus we have covered all of the subalgebras and our proof is complete.
\end{proof}
\end{theorem}


\newpage
\begin{thebibliography}{99}
\bibitem{Dixon} Dixon,G.M. Division Algebras: Octonions, Quaternions,
Complex Numbers and the Algebraic Design of Physics. Kluwer Academic Publishers,1994
\bibitem{Ebbinghaus88} Ebbinghaus,H.-D. et al. Numbers. Springer-Verlag, 1988
\bibitem{Freudenthal65} Freudenthal,H. Lie Groups in the Foundations of
Geometry, Adv. Math I 135 (1965).
\bibitem{Gilmore} Gilmore,R. Lie Groups, Lie algebras and Some of Their
Applications, John Wiley and Sons, 1941.
\bibitem{GurseyTze96} G\"{u}rsey,F. and Tze,C. On the Role of
Division, Jordan and Related algebras in Particle Physics, World
Scientific, 1996
\bibitem{Howe} Howe,R. and Umeda,T. The Capelli identity, the double
commutant theorem, and multiplicity-free actions. Math. Ann. 290 (1991).
\bibitem{JacobsonELA} Jacobson,N. Exceptional Lie Algebras, Marcel Dekker, Inc. 1971
\bibitem{Kantorsolod} Kantor,I.L and Solodovnikov,A.S. Hypercomplex
numbers; An elementary introduction to Algebras, Springer-Verlag, 1980
\bibitem{Molev} Molev,A. and Nazarov,M. Capelli Identities for Classical
Lie Algebras, Math. Ann. 313 (1999).
\bibitem{Nazarov} Nazarov,M. Capelli Elements in the classical universal
enveloping algebras, math/9811129.
\bibitem{Porteous} Porteous,I.R. Topological Geometry, Van Nostrand
Reinhold Company Ltd, 1969.
\bibitem{Ramond} Ramond,P. Introduction to Exceptional Lie  Groups and
Algebras, Preprint, CALT-68-577, 1976
\bibitem{Schafer66} Schafer,R.D. Introduction to Non-associative
Algebras, Academic Press, 1966
\bibitem{Sudbery84} Sudbery,A. Division Algebras, (Pseudo-) Orthgonal
Groups and Spinors, J.Phys A 17 (1984).
\bibitem{Tits62} Tits,J. Algebres Alternatives, Algebres de Jordan et
Algebres de Lie Exceptionelles, Proc. Colloq Utrecht 135 (1962).
\bibitem{Weyl} Weyl,H. The Classical Groups, Princeton University Press, 1946
\end{thebibliography}
\end{document}